\documentclass[11pt,a4paper]{article}
\usepackage{caption2}
\usepackage{amssymb}
\usepackage{color}
\usepackage{amsmath,amsthm}
\usepackage{graphicx}
\usepackage{empheq}
\usepackage{indentfirst}
\usepackage{hyperref}
\usepackage{float}

 \newtheorem{thm}{Theorem}[section]
 \newtheorem{cor}{Corollary}[section]
 \newtheorem{lem}{Lemma}[section]
 
 \theoremstyle{definition}
 \newtheorem{defn}{Definition}[section]
 \newtheorem{rem}{Remark}[section]
 \newtheorem{example}{Example}[section]
 \numberwithin{equation}{section}
\def\f{\frac}
\def\pa{\partial}

\def\al{\alpha}
\def\e{\eqref}

\def\lab{\label}

\def\i1n{i=1,\cdots,n}
\def\j1n{j=1,\cdots,n}
\def\ij1n{i,j=1,\cdots,n}

\def\R{\mathbb R}

\newcommand{\be}{\begin{equation}}
\newcommand{\ee}{\end{equation}}
\newcommand{\beq}{\begin{equation*}}
\newcommand{\eeq}{\end{equation*}}

\begin{document}
\title{A one-step optimal energy decay formula for indirectly  nonlinearly damped hyperbolic systems coupled by velocities}
\author{
Fatiha Alabau-Boussouira\thanks{IECL, Universit\'{e} de Lorraine and CRNS (UMR 7502),
57045 Metz, France. E-mail: {\tt fatiha.alabau@univ-lorraine.fr}.  
},
 Zhiqiang Wang\thanks{School of Mathematical Sciences and Shanghai Key Laboratory for Contemporary
Applied Mathematics, Fudan University, Shanghai 200433, China. E-mail: {\tt wzq@fudan.edu.cn}. 
}
\  and
Lixin Yu\thanks{
School of Mathematics and Information Sciencs, Yantai University,
Yantai 264005, China. E-mail: {\tt fdylx01@sina.com.cn}.
} 
}
%


\maketitle


\begin{abstract}
In this paper, we consider the energy decay of  a damped hyperbolic system of wave-wave type 
which is coupled through the velocities. We are interested in the asymptotic properties of the solutions of this system
in the case of indirect nonlinear damping, i.e. 
when only one equation is directly damped by  a nonlinear  damping. 
We prove  that the total energy of the whole system decays as fast as the damped single equation.
Moreover, we give a one-step general explicit decay formula for arbitrary nonlinearity. 
Our results shows that the damping properties are fully transferred from the damped equation to the
undamped one by the coupling in velocities, different from the case of couplings through displacements as shown in \cite{AB01, ACK01, AB02, AL12} 
for the linear damping case, and in \cite{AB07} for the nonlinear damping case. 
The proofs of our results  are based on multiplier techniques, weighted nonlinear integral inequalities and the optimal-weight convexity method of \cite{AB05, AB10}.
\end{abstract}



 \section{Introduction }
  Let $\Omega$ be a bounded subset of $\mathbb{R}^n$ with a smooth boundary
denoted by $\Gamma$.  We consider  the following wave system
 \be
 \lab{wave}
  \left\{
 \begin{array}{l}
  u''-\Delta u+\alpha(x)v'+\rho (x,u') =0 \quad\quad\  \text{in } \ \Omega\times (0,+\infty),\\
  v''-\Delta v-\alpha(x) u'=0\quad\quad\quad \text{in } \ \Omega\times (0,+\infty),\\
  u=v=0\quad\quad\quad\quad\text{on} \  \Gamma\times  (0,+\infty),\\
  (u,u')(0)=(u^0,u^1),  \  ( v,v')(0)=(v^0,v^1)\quad\quad\text{in} \  \Omega\,.
  \end{array} 
 \right.
 \ee 
 To a strong solution of this system, we associate
the energy defined by
  \be\lab{ene-wave} 
  E(t)=\f12\int_\Omega (|u'|^2+|\nabla u|^2+|v'|^2+|\nabla v|^2)dx \,.
 \ee
 One can show that the energy of the strong solutions of this system satisfies
 \be\lab{ene-wave-dissip} 
E'(t) =-\int_{\Omega} u' \rho(x, u')dx \,, 
\ee

%
In the sequel, the assumption on $\rho$ will ensure that $\rho(\cdot, s)s  \geqslant  0$ a.e. in $\Omega$ and for all $s \in \mathbb{R}$, so that
the energies of the strong solutions satisfy $E'(t)\leq0$. Hence the nonlinear term in our coupled system is
indeed a damping term, so that one expects the energies to decay to $0$ at infinity.

 Let us recall that for the scalar damped wave equation, that is for
 $$
 u''-\Delta u+\rho (x,u') =0 \quad \text{in } \ \Omega\times (0,+\infty)\,, u=0 \quad   \text{on} \  \Gamma\times  (0,+\infty)\,.
 $$
 It is well known that when the damping term is linear, i.e. when $\rho(.,u')=a(\cdot)u'$, where $a  \geqslant  0$ a.e. on $\Omega$, the energy of
  the solution decays exponentially under some  some geometric  conditions on the support of $a$ (see \cite{BLR92, Lio}).
When the damping term $\rho(\cdot, \cdot)$ is nonlinear with respect to the second variable, 
  \cite{AB05} and \cite{AB10}, \cite{ABAM11}
  give, under some suitable geometric conditions, a one-step explicit energy decay formula in terms of the behavior of the nonlinear feedback
  close to the origin. These results rely on a general weighted nonlinear  integral
inequality together with an optimal-weight convexity method developed in \cite{AB05}.
 If no geometric assumptions on the damping region are made, the decay is known to be of logarithmic type for a linear damping (see e.g. \cite{Bell,Fu11}).
 
Let us go back to the above wave-wave coupled system \eqref{wave}. At this stage, four main features characterize this system: 
\begin{itemize}
\item only the first equation is damped
\item the damping $\rho$ may have an arbitrary growth around 0 with respect to the second variable 
\item the coupling coefficient $\alpha$ may vanish in some parts of $\Omega$
\item the coupling is acting through the velocities.
\end{itemize}
Let us now comment  on these features.

$\bullet$ The fact that only one equation of the coupled system is damped refers to the so-called class of "indirect" stabilization problems initiated and studied in \cite{AB01, ACK01, AB02} and further studied in \cite{AB06, ACG11, AL12}. Indeed, when dealing with coupled systems, it may be impossible or too expensive to damp each equation. Such an example is provided for instance by the Timoshenko system \cite{AB07,kim, soufyane}. More generally, coupled system involving some undamped equations, are said to be indirectly damped. From the point of view of applications in control theory, a challenging question is to determine whether the single feedback is sufficient to guarantee that the energy of the full system decays to $0$ at infinity and to determine at which rate. In this latter case, the lack of feedback on the second equation is compensated by the coupling effects.

$\bullet$ The case of general damping feedbacks, that is with arbitrary growth close to $0$, has received a lot of attention since more than a decade. The first result in this direction has been derived in \cite{LT93}, however no general simple explicit formula was provided except for linearly or polynomially growing dampings close to $0$. Such first examples of explicit general formula are given in \cite{Mart} (see also \cite{liuzuazua}), but this formula does not allow to recover in a single step the expected quasi-optimal energy decay rate in the polynomial case (or for polynomial-logarithmic growth). As far as we know, the first result giving a general one-step quasi-optimal semi-explicit formula is given in \cite{AB05}. A further analysis based on a suitable and original classification of the feedback growth has been introduced in \cite{AB10}. This classification gives a very simple one-step explicit energy decay formula for general feedbacks growth, provided that this growth is not close to a linear behavior. A more complex semi-explicit formula holds in the general case including feedbacks  with a growth close to a linear one around $0$. By a {\em one-step} formula, we mean here, a formula 
which gives a decay rate depending explicitly on the feedback in a simple explicit way. 
In particular, this formula does not require further steps as in most of the existing literature to lead to explicit expressions. Moreover the optimality of the formula is proved for the corresponding finite dimensional systems or for semi-discretized scalar wave or plate equations, whereas optimality results are proved for some examples in the infinite dimensional case in \cite{AB05}, using results of \cite{VCM, VCM2}. For previous results
and also results for semilinear damped equations, one can see \cite{chen, kom,  nakao, zuazua, zuazua2}.

$\bullet$ When the coupling coefficient is bounded below by a strictly positive constant, the coupling is active in the whole domain $\Omega$, so that the equations are coupled in the whole domain. If $\alpha$ is nonnegative on $\Omega$ but is allowed to vanish on a subset of $\Omega$, the equations are "uncoupled" in $\Omega \backslash supp\{\alpha\}$, so that  in this region the second equation is decoupled from the first equation and is undamped. Such cases are harder to handle. A first study in this framework, but for linearly damped wave-type equations coupled in {\bf displacements} is given in \cite{AL12}.

$\bullet$ When indirect damping occurs through displacements, that is for systems coupled in displacements, as for instance 
$$\left\{\begin{array}{l}
u''-\Delta u + \alpha(x)v+a(x)u'=0\qquad\ \
\text{in}\  (0,\infty)\times \Omega,\\
v''-\Delta v +  \alpha(x)u=0\qquad  \text{in}\  (0,\infty)\times \Omega,\\
u=v=0  \qquad \text{on}\  (0,\infty)\times \Gamma,\\
(u(0), u'(0) )=(u^0, u^1), (v(0),v'(0))=(v^0,v^1),\qquad  \text{in}\  \Omega,\end{array}\right.$$ 
it has been shown in \cite{ACK01} that even for constant coefficients $\alpha$, the energy of this
linearly damped system never decays exponentially, but decays only polynomially with a  decay rate depending on the smoothness of the initial data and a general lemma announced in \cite{AB01} (see \cite{ ACK01,AB02} for a proof). These results are based on the method of higher order energies initiated in \cite{AB01} and developped in \cite{ACK01} and \cite{AB02} for the indirect boundary damping cases in an abstract setting and applied to various examples.
 The result of \cite{ACK01} has been generalized to coefficients $\alpha$ that vanishes on a subset of $\Omega$
in \cite{AL12}, under certain assumptions on the supports of $\alpha$ and
$a$ (roughly speaking they are both supposed to satisfy the Geometric Control Condition (GCC) of \cite{BLR92}). Further results on coupled models with distributed dampings but satisfying hybrid boundary conditions have been obtained in \cite{ACG11}. Sharper results have been obtained through interpolation techniques extending the first results of \cite{AB01}.
Let us further mention that in \cite{Fu}, the author shows that the energy decays logarithmically  under the assumption 
$supp\{\alpha\}\cap supp\{a\}\neq \emptyset$. 
Hence when the coupling acts through displacements, indirect stabilization occurs but in a weaker form than the one of the corresponding scalar case, since exponential stabilization does not hold even for a linearly damped case.

The goal of this paper is to generalize the quasi-optimal energy decay formula given for {\bf scalar} wave-type systems in \cite{AB10} (see also \cite{AB05}), to the case of  {\bf coupled} systems in velocities, under the above four features. More precisely, we prove, under some geometric conditions on 
 the localized damping domain and the localized coupling domain, 
 that  the energy of these kinds of  system decays as fast as that of the corresponding scalar nonlinearly damped equation. Hence, the coupling through 
 velocities allows a full transmission of the damping effects, quite different from the coupling through the displacements.

%

 The optimality of the above estimates has been proved for finite dimensional equations, including the semi-discretized wave equations in \cite{AB10}. In the infinite dimensional setting, lower energy estimates or optimality are open questions. Optimality has been only proved in the particular case of one-dimensional wave equation with boundary damping (see \cite{AB05, VCM, VCM2}).  Lower energy estimates have been established in \cite{AB10, AB10bis, AB11} for scalar one-dimensional wave equations, scalar Petrowsky equations in two-dimensions and Timoshenko systems. We use the comparison method developed in \cite{AB11} to extend these results to one-dimensional wave systems coupled by velocities. 
 
 \begin{rem}
The method presented here is general and can easily be adapted to handle corresponding coupled systems of plate equations, elasticity models,
and more complex examples in the spirit of the general approach given in \cite{AB05}. Our aim through this paper is to give a general methodology on a concrete PDE example to show that if the damping effects are suitably transferred through the coupling operators, then indirect stabilization can produce damping mechanisms of the quality of a direct damping for the corresponding scalar equation.

Note that one can also present, with no additional mathematical originality and no gain with respect to applications, all our results by means of a more "Lyapunov" presentation. In this case,
it is sufficient not to integrate over the time interval and to handle terms of the form $\dfrac{d}{dt}\big[\int_{\Omega} u'^2dx\big]$ for instance, and multiply afterwards by a weight function, which can be a weaker (and less good weight function) than in the original method introduced for the first time in \cite{AB05} and announced in \cite{AB04}. 
\end{rem}

 The rest of this paper is organized as follows:
In Section \ref{sec2}, we give some basic preliminaries, assumptions and notations. 
The main results, including Theorems \ref{thm-wellposed} on well-posedness of \eqref{wave} ,
Theorem \ref{thm-wave} on energy decay for polynomially growing damping case,
Theorem  \ref{thm-wave-gen_decay} on energy decay  for general nonlinear damping case, 
Theorem \ref{thm-wave-lower}  on lower  energy estimates for one dimensional system \eqref{wave1D}, are presented in Section \ref{sec3}.
Explicit decay rates corresponding to some typical dampings are also provided in Subsection \ref{exdecayrate}. 
As the main tool of deriving the quasi-optimal one-step explicit energy decay formula,
 the optimal-weight convexity method together with general weighted nonlinear integral inequalities are introduced in Section \ref{sec4}.
The proof of  the main results, as well as the decay rates of Example \ref{ex3.1}-\ref{ex3.4} are given in Section \ref{sec5}.
 Finally, Section \ref{sec6} is used to prove Lemma \ref{lem-wave} on weighted energy estimates  for a single non homogeneous wave equation by the multiplier method.

\section{Preliminaries, assumptions and notations} \lab{sec2}


\begin{description}

\item {\small {\bf Notations}}

For the simplicity of statement, we denote in the whole paper  $  L^2(\Omega) $ by $L^2$,  $ H_0^1(\Omega)  $ by $H^1_0$, $  H^2(\Omega)  $ by $H^2$. 
Moreover, we say that the initial data are in the energy space whenever $(u^0,u^1)\in  H_0^1 \times   L^2 $ and 
$(v^0,v^1)\in  H_0^1 \times   L^2 $ and the  initial data are smooth if  $(u^0,u^1)\in
  (H^2  \cap  H_0^1) \times  H_0^1 $ and $(v^0,v^1)\in   (H^2  \cap  H_0^1) \times  H_0^1 $.

\item {\small {\bf Geometric conditions}}


As already mentioned in the introduction, stabilization results for wave-like systems require geometric conditions on the region where the feedback is active. In the sequel, we shall consider the so-called Piecewise Multipliers Geometric Condition (denoted by PMGC, in short):

\begin{defn}[\bf PMGC] We say that a subset $\omega\subset\Omega$ satisfies
the PMGC, if there exist subsets $\Omega_j\subset\Omega$ having Lipschitz
boundaries and points $x_j\in \mathbb{R}^N,\ j=1,\cdots,J$, such that
$\Omega_i\cap\Omega_j=\emptyset$ for $i\neq j$ and $\omega$ contains a
neighborhood in $\Omega$ of the set
$\cup_{j=1}^J\gamma_j(x_j)\cup(\Omega\setminus
\cup_{j=1}^J\Omega_j)$, where $\gamma_j(x_j)=\{x\in
\partial\Omega_j: (x-x_j)\cdot \nu_j(x) \geqslant  0\}$ and $\nu_j$
is the outward unit normal to $\partial\Omega_j.$
\end{defn}

\item {\small  {\bf Assumptions on the coupling coefficient}}
 
 \medskip

We assume that the coupling function $\alpha \in C(\bar{\Omega})$   satisfies
 \be \lab{HC}
 \tag{HC}
 \left\{
 \begin{array}{l}
 \exists\  \alpha_+>0,   \alpha_->0 \ \text{such that} \\
 \alpha_+  \geqslant  \alpha(x) \geqslant   0,\quad  \forall x\in  \Omega,\\
  \alpha(x) \geqslant  \al_->0,\quad  \forall x\in \omega_c \subset \Omega,
 \end{array}
 \right.
 \ee 
where 
and $\omega_c$ is an open subset of $\Omega$ with positive measure.

\item {\small {\bf Assumptions on the feedback}}

We consider feedbacks $\rho$ with an arbitrary growth close to $0$. However, to give the reader a better insight of the scope and challenge
of one-step explicit general quasi-optimal energy decay formulas, we first provide the result and proof for polynomially growing feedbacks, for which
the proofs are easier and then the general result for arbitrary growing feedbacks. Hence, we detail below the two sets of assumptions: the one for the polynomial case, then those for the general case.

The assumptions in the case of {\bf polynomially} growing feedbacks is as follows 
 \be \label{HF}
 \tag{HF$_p$} \quad
 \left\{\begin{array}{l}  
 \rho\in C(\overline{\Omega} \times \mathbb{R}) \,, \rho(x, 0)= 0 \quad \forall x\in \Omega\,,\\
 s \mapsto  \rho(x,s) \mbox{ is nondecreasing } \forall \ x \in \Omega,\\
\exists \  c >0 \ \text{and}\   p \geqslant  1\ , \exists \ a \in C(\overline{\Omega}) \text{ such that} \\
a(x)|s|\leqslant |\rho(x, s)|\leqslant c a(x) |s|,\quad\forall x \in \Omega\,, |s| \geqslant  1,\\
a(x) |s|^p\leqslant |\rho(x, s)|\leqslant c a(x) |s|^{\f{1}{p}} ,\quad\forall x \in \Omega\,, |s|\leqslant 1, \text{ where }\\
a  \geqslant  0 \text{ on } \Omega \,,
\exists \  a_->0 \text{ such that }  \ a(x) \geqslant  a_-,\quad\forall  x\in \omega_d \subset \Omega,\end{array}\right.
 \ee
where $\omega_d$ is an open subset of $\Omega$ with positive measure.

The assumptions in the case of {\bf arbitrary} growing feedbacks is as follows
%
\be \label{HFg}
 \tag{HF$_g$} \quad
 \left\{\begin{array}{l}   
 \rho\in C(\overline{\Omega} \times \mathbb{R}) \,, \rho(x, 0)= 0 \quad \forall x\in \Omega\,,\\
 s \mapsto  \rho(x,s) \mbox{ is nondecreasing } \forall \ x \in \Omega,\\
\exists \  c >0,   \exists \ a \in C(\overline{\Omega}) \mbox{ and } \exists \ g \in \mathcal{C}^1(\mathbb{R}) \text{ such that }\\
a(x)|s|\leqslant |\rho(x,s)|\leqslant ca(x) |s|,\quad\forall x \in \Omega\,, |s| \geqslant  1,\\
a(x)g(|s|)\leqslant |\rho(x,s)|\leqslant c a(x)g^{-1}(|s|) ,\quad\forall x \in \Omega\,, |s|\leqslant 1, \text{ where }\\
a  \geqslant  0 \text{ on } \Omega \,,
\exists \ a_->0 \text{ such that } a(x) \geqslant  a_-,\quad\forall  x\in \omega_d \subset \Omega,\\
 g \mbox{ is a strictly increasing and odd function}.
\end{array}\right.
 \ee
\begin{rem}
Thanks to the hypotheses \eqref{HF} or \eqref{HFg}, we have
\begin{equation}\label{damprho}
\rho(x,s)s  \geqslant   0,\quad \forall x\in \Omega\,, \forall \ s \in \mathbb{R}\,,
\end{equation}
\noindent which ensures that the energy of the solutions of the above wave system is nonincreasing.
\end{rem}
\begin{rem}
 Note that  we can infer from $(HF_g)$ that for very $\varepsilon \in (0,1)$, there exists constants $c_1>0, c_2>0$ such that
 \begin{equation}\label{gwrho}
  \left\{\begin{array}{l}   
c_1a(x)|s|\leqslant |\rho(x,s)|\leqslant c_2 a(x)|s|,\quad\forall x \in \Omega\,, |s| \geqslant  \varepsilon,\\
c_1a(x)g(|s|)\leqslant |\rho(x,s)|\leqslant c_2 a(x)g^{-1}(|s|) ,\quad\forall x \in \Omega\,, |s|\leqslant \varepsilon\,.
\end{array}\right.
 \end{equation}
 \end{rem}

 {\bf Convexity assumptions on the feedback and some definitions}
 
 Assume that  \eqref{HFg} holds. Then, following \cite{AB05, AB10}, we assume that the function $H$ defined by
\begin{equation}\label{defH}
H(x)=\sqrt{x}g(\sqrt{x}),
\end{equation}
is strictly convex in a right neighborhood of $0$, i.e. on $[0,r_0^2]$ for some sufficiently small $r_0\in(0,1]$. 
We define the function $\widehat H$ on $\R$ by $\widehat H(x)=H(x)$ for every $x\in[0,r_0^2]$ and by $\widehat H(x)=+\infty$ otherwise, and we define the function $L$ on $[0,+\infty)$ by
\begin{equation}\label{defL}
L(y)=\left\{\begin{array}{ll}
\displaystyle{\frac{\widehat H^\star(y)}{y}} & \textrm{if }y>0,\\
0 & \textrm{if }y=0,
\end{array}
\right.
\end{equation}
where $\widehat H^\star$ is the convex conjugate function of $\widehat H$, defined by
$$
\widehat H^\star (y)=\sup_{x\in \R}\{xy-\widehat H(x)\}.
$$
By construction, the function $L:[0,+\infty)\rightarrow[0,r_0^2)$ is continuous, one-to-one, onto and increasing, moreover it is easy to check that
\begin{equation}\label{LH'}
0 < L(H'(r_0^2)) <r_0^2 \,
\end{equation}
holds (see \cite{AB05, AB10} for a complete proof).
We also define the function $\Lambda_H$ on $(0,r_0^2]$ by
\begin{equation}\label{defLambdaH}
\Lambda_H(x)=\frac{H(x)}{xH'(x)}.
\end{equation}
\begin{rem}
The function $\Lambda_H$ has been introduced for the first time by the first author in \cite{AB10}. It is an essential tool
to classify the feedback growths around $0$ and to simplify the decay estimate formula given in \cite{AB05} --without loosing
optimality properties-- for the feedbacks having a growth around $0$ which is not close to a linear one (as explained below). 
\end{rem}
\begin{rem}
Note that due to our convexity assumptions, we have
$$
\Lambda_H([0,1]) \subset [0,1]\,.
$$
\noindent When $g(x)=x$ (linear feedback case), $\Lambda_H(.) \equiv 1$. If we now set for instance, $g(x)=x (\ln(1/x))^{-q}$ in $(0, \varepsilon]$ 
with $q>0$ and $\varepsilon>0$, then $\limsup \limits_{x \rightarrow 0^+} \Lambda_H(x)=1$. Many other examples of feedbacks such that 
$\limsup \limits_{x \rightarrow 0^+} \Lambda_H(x)=1$ can be given, they characterize a growth which is close to a linear one. This leads to the following definition:
\begin{defn}
We say that a feedback $\rho$ satisfying \eqref{HFg} has a growth close to a linear one in a neighborhood of $0$, if it is such that the function $H$ defined by \eqref{defH} satisfies  $\limsup \limits_{x \rightarrow 0^+} \Lambda_H(x)=1$. Otherwise, one says that the feedback $\rho$ is away from a linear growth. 
\end{defn}
On the opposite side, for functions $g$ which converge very fast to $0$ as $x$ goes to $0$, such as for instance $g(x)=e^{-1/x}$ for $x \in (0,\varepsilon]$
(and many other examples), one has $\limsup\limits_{x \rightarrow 0^+} \Lambda_H(x)=0$. 

For polynomially growing feedbacks, e.g. when $g(x)=x^p$ with $p>1$, we have
$\Lambda_H(.) \equiv \frac{2}{p+1}$. For feedbacks such as $g(x)=x^p (\ln{1/x})^q$ with $p>1\,,q>0$, we still
have $\limsup\limits _{x \rightarrow 0^+} \Lambda_H(x)=\frac{2}{p+1}$.

We will see later on, in Theorem~\ref{thm-wave-gen_decay}, that the case of feedbacks close to a linear behavior as $x$ goes to $0$ has to be distinguished from the other cases. 
\end{rem}
Finally, we define, for $x\geq 1/H'(r_0^2)$,
\begin{equation}\label{defpsir}
\psi_0(x)=\frac{1}{H'(r_0^2)}+\int_{1/x}^{H'(r_0^2)}\frac{1}{\theta^2(1-\Lambda_H((H')^{-1}(\theta)))}\,d\theta.
\end{equation}
\begin{rem}
Note that when $g'(0) \neq 0$,
$g$ has a linear growth close to $0$. Therefore, this case is similar to the linear case which is already well-known. We thus focus in the sequel on
the cases where $g'(0)=0$.
\end{rem}

\end{description}

\section{Main results} \label{sec3}

\subsection{Well-posedness}

We set $\mathcal{H}=( H_0^1 \times   L^2 )^2$ and set $U=(u,p, v,q)$. We equip $\mathcal{H}$ with the scalar product
$$
\langle U \,, \widetilde{U} \rangle = \int_{\Omega} \Big(\nabla u \cdot \nabla  \widetilde{u} + p  \widetilde{p} +
\nabla v \cdot \nabla  \widetilde{v} + q  \widetilde{q}  \Big) dx
$$
\begin{thm} \label{thm-wellposed}
Assume \eqref{HFg} and that $\alpha \in L^{\infty}(\Omega)$. Then for all initial data in energy space, there exists a unique solution 
$(u,v) \in \mathcal{C}([0,+\infty);( H_0^1 )^2)\cap\mathcal{C}^1([0,+\infty);(  L^2 )^2)$
of \eqref{wave}.
Moreover for any smooth initial data, the solution satisfies
$(u,v) \in  L^{\infty}([0,+\infty);(  H^2  \cap  H_0^1 )^2)\cap
W^{1, \infty}([0,+\infty);( H_0^1 )^2) \cap W^{2, \infty}([0,+\infty);(  L^2 )^2)$.
Moreover, in this latter case the energy of order one, defined by
\be\lab{ene-wave-1} 
  E_1(t)=\f12\int_\Omega (|u_{tt}|^2+|\nabla u_t|^2+|v_{tt}|^2+|\nabla v_t|^2)dx \,,
 \ee
is non increasing, i.e., 
\be\lab{ene-wave-1-dissip} 
E_1(t)\leq E_1(0)\,.
\ee
\end{thm}

\subsection{One-step quasi-optimal energy decay rate for the wave-wave system}

{\bf First case:  polynomiaily growing dampings close 0}

For the sake of clarity, we first provide the results in the case of a polynomially growing feedback $s \mapsto \rho(\cdot, s)$.

\begin{thm} \label{thm-wave}
Assume that $p>1$, \eqref{HF} and \eqref{HC} hold. Assume also that
$\omega_d$ and $\omega_c$ satisfy the PMGC.  Then there exists $\alpha^{\ast}>0$ such that for any $\alpha_+ \in (0,\alpha^{\ast}]$ and
any non vanishing initial data in the energy space, the total energy of \eqref{wave}
defined in \eqref{ene-wave} 
decays as 
 \be \lab{decay-poly} 
  E(t)\leqslant C_{E(0)} \, t^{\f{2}{1-p}}, \quad  \forall t\in [T_{E(0)},+\infty), 
   \ee
where $C_{E(0)},T_{E(0)}>0$  are constants  depending on $E(0)$.
 \end{thm}
 
\begin{rem}
If $p=1$,  we can follow the proof of Theorem \ref{thm-wave} 
and  obtain, by using  Lemma \ref{lem-expo} instead of  Lemma \ref{lem-poly},
 the exponential stability for system \eqref{wave} 
 \be \lab{decay-expo} 
  E(t)\leqslant 
  C E(0) \,  e^{-\kappa t}, \quad \forall t\in [0,+\infty),  \ee
where  $C>0, \kappa>0$  are constants  independent of the initial data.
\end{rem}

{\bf Second case: arbitrarily growing dampings close to the origin}

\medskip

If we consider a more general nonlinear damping $\rho$, we provide below a quasi-optimal one-step explicit energy decay formula following
the optimal-weight convexity method together with general weighted nonlinear integral inequalities developed in \cite{AB05, AB10}.

\begin{thm} \label{thm-wave-gen_decay}
Assume \eqref{HFg} and \eqref{HC} hold. Assume that the function $H$ defined by \eqref{defH} is strictly convex in $[0,r_0^2]$ for some sufficiently small $r_0\in(0,1]$ and satisfies $H(0)=H^{\prime}(0)=0$. We define the maps $L$ and $\Lambda_H$ respectively as in \eqref{defL} and \eqref{defLambdaH}.
Assume also that $\omega_d$ and $\omega_c$ satisfy the PMGC. Then there exists $\alpha^{\ast}>0$ such that for any $\alpha_+ \in (0,\alpha^{\ast}]$ and
any non vanishing  initial data  in the energy space, the total energy of  \eqref{wave}
defined in \eqref{ene-wave} 
decays as
\begin{equation}\label{decayenergymain}
\displaystyle{E(t)\leq2 \beta_{E(0)} L\Big(\frac{1}{ \psi_{0}^{-1}(\frac{t}{M})}\Big) \ , \quad \forall \ t\geq
\frac{M}{H'(r_0^2)}}\,,
\end{equation}
\noindent where $\beta_{E(0)}$ is defined by \eqref{defbeta}, $M$ is defined
by \eqref{M} and independent of $E(0)$.
Furthermore, if $\limsup \limits_{x \rightarrow 0^+} \Lambda_H(x)<1$, then $E$ satisfies the following simplified decay rate
\begin{equation}\label{decayenergysimpmain}
\displaystyle{E(t)\leq2 \beta_{E(0)} \Big(H'\Big)^{-1}\Big(\frac{\kappa M}{t}\Big) \ , }
\end{equation}
\noindent for $t$ sufficiently large, and where $\kappa>0$ is a  constant independent  of $E(0)$. 
 \end{thm}
 \begin{rem}
 Note that when $g(x)=x^{1/p}$ wih $p>1$, then $H(x)=x^{(p+1)/2}$ so that $\limsup \limits_{x \rightarrow 0^+} \Lambda_H(x)=\frac{2}{p+1}<1$. Then
 the formula \eqref{decayenergysimpmain} gives back the energy decay rate of $t^{-2/(p-1)}$ given in Theorem~\ref{thm-wave}. 
 \end{rem}

  \begin{rem}
 The smallness of $\alpha$ can be reduced if we assume additionally $supp \{\alpha\} \subset \omega_d$ in Theorem \ref{thm-wave-gen_decay}. 
 Actually, we can choose $\delta=1$ in the proof of \e{newb3}, and in that case, the last term of \e{newb3} can be replaced by 
$$   C_1 \int_S^T\phi(t)\int_{\omega_d}  |u' |^2dx dt  $$
Then \e{decay-dominant-energ}  follows easily by choosing $\varepsilon_1>0$ sufficiently small in \e{newb3}. The proof afterwards is the same.
 \end{rem}
 
 \subsection{Lower energy estimates}
 
The optimality of the above estimates are open questions. 
Here we use the comparison method developed in \cite{AB11} to establish the lower estimate of the energy of the one-dimensional coupled wave system.
Actually, we consider the following wave-wave system in $\Omega=(0,1) \subset \mathbb{R}$
and $\rho(x,s)=a(x)g(s)$ for all $ x \in \Omega$ and all $s \in \mathbb{R}$.   
 \be
 \lab{wave1D}
  \left\{
 \begin{array}{l}
  u_{tt}-u_{xx}+\alpha(x)v_t+a(x)g(u_t) =0 \,,\quad\  \text{in } \ \Omega\times (0,+\infty),\\
  v_{tt}- v_{xx}-\alpha(x) u_t=0  \,, \quad \text{in } \ \Omega\times (0,+\infty),\\
  u(t,0)=u(t,1)=v(t,0)=v(t,1)=0  \,, \quad \text{for } \   t \in (0,+\infty),\\
  (u,u_t)(0,x)=(u^0(x),u^1(x)),  \  ( v,v_t)(0,x)=(v^0(x),v^1(x)) \,,\quad\text{for } \  x \in \Omega\,.
  \end{array} 
 \right.
 \ee 
 We define $H$ by \eqref{defH}, $\Lambda_H$ by \eqref{defLambdaH} and  consider the following assumptions
\be \label{HFl}
\tag{HF$_l$} \begin{cases}
\exists r_0>0 \mbox{ such that the function } H: [0,r_0^2] \mapsto \mathbb{R}
\mbox{ defined by } \eqref{defH}\\
\mbox{ is strictly convex on } [0,r_0^2]\,  \mbox{ and }\\
 \mbox{ either } \ 
0< \liminf \limits_{x \rightarrow 0+}\Lambda_H(x)
\le \limsup \limits_{x \rightarrow 0+}\Lambda_H(x) <1 \,\\
 \mbox{ or there exist } \mu>0   \mbox{ and } z_1
\in (0,z_0] \, \mbox{with} \ \eqref{defz0}  
\mbox{ such that } \\
\qquad 0< \liminf \limits_{x \rightarrow 0+}
\Big(\frac{H(\mu\, x)}{\mu \, x}\int_x^
{z_1}\frac{1}{H(y)}\,dy
\Big)\,
\mbox{ and } 
\limsup \limits_{x \rightarrow 0+}\Lambda_H(x) <1.
\end{cases}
\ee
\noindent

 \begin{thm} \label{thm-wave-lower}
Assume that \eqref{HFg}, \eqref{HFl}, and \eqref{HC} hold. Assume also that
$\omega_d$ and $\omega_c$ satisfy the PMGC, and that $\alpha_+ \in (0,\alpha^{\ast}]$ where $\alpha^{\ast}>0$ is as in
Theorem \ref{thm-wave-gen_decay}. Then for all non vanishing smooth initial data,
 there exist $T_0>0$ and $T_1>0$ such that the energy of \eqref{wave1D} satisfies the lower estimate
\begin{equation}\label{carac-low-wave}
  E(t)  \geq \frac{1}{\gamma_{s}^2 C_{s}^2} \Big(\big(H^{\prime}\big)^{-1}\Big(\frac{1}{t-T_0}\big)\Big)^2\,,  \quad \forall \ t \ge T_1+T_0 \,, 
\end{equation}
\noindent where $\gamma_{s}=4 \sqrt{E_1(0)}$, $C_s$ is as in Lemma \ref{compar}.
 \end{thm}
 
 \subsection{Some examples of decay rates}\label{exdecayrate}

For the sake of completeness,  we give some significative examples (taken from \cite{AB05, AB10}) of feedback growths together with the
resulting energy decay rate when applying our results. 
In the sequel, $C_{E(0)}>0$ stands for a constant depending on  $E(0)$,
while $C^{\prime}_{E_1(0)}>0$ is a constant depending  on  $E_1(0)$.
\begin{example}\label{ex3.1}[The Polynomial Case]
Let  $g\left( x\right) =x^{p}$   on $(0,r_{0}]$ with $p>1$.
Then the energy of \eqref{wave} satisfies the
estimate 
\begin{equation*}
E\left( t\right)\leq C_{E\left( 0\right) }t^{-\frac{2}{p-1}},
\end{equation*}%
for $t$ sufficiently large and for all non vanishing initial data in the energy space. 
Moreover,   the energy of  \eqref{wave1D} satisfies the
estimate 
\begin{equation*}
E\left( t\right)  \geqslant   C^{\prime}_{E_1(0)} t^{-\frac{4}{p-1}},
\end{equation*}%
for $t$ sufficiently large and for all non vanishing smooth initial data.
\end{example}

\begin{example}\label{ex3.2}[Exponential growth of the feedback]
Let $g\left( x\right) =e^{-\frac{1}{x^2}}$ on $(0,r_{0}].$ Then
the energy of  \eqref{wave}  decays as
\begin{equation*}
E\left( t\right)\leq C_{E\left( 0\right) }\left( \ln \left( t\right)
\right) ^{-1},
\end{equation*}%
for $t$ sufficiently large and for all non vanishing initial data in the energy space. 
Moreover, the energy of \eqref{wave1D} satisfies the
estimate 
\begin{equation*}
E\left( t\right)  \geqslant   C^{\prime}_{E_1(0)}  \left( \ln \left( t\right)
\right) ^{-2},
\end{equation*}%
for $t$ sufficiently large and for all non vanishing smooth initial data.

\end{example}

\begin{example}\label{ex3.3}[Polynomial-logarithmic growth]
Let $g\left( x\right) =x^{p}\left( \ln \left( 1/x\right)
\right) ^{q}$ on $(0,r_{0}]$ with $p>2$ and $q>1$.  Then the energy of \eqref{wave} decays as 
\begin{equation*}
E\left( t\right)\leq C_{E\left( 0\right) } t^{-2/(p-1)} (\ln (t))^{-2q/(p-1)}
\end{equation*}%
for $t$ sufficiently large and for all non vanishing initial data in the energy space. 
Moreover the energy of  \eqref{wave1D} satisfies the
estimate 
\begin{equation*}
E\left( t\right)  \geqslant  C^{\prime}_{E_1(0)}  t^{-4/(p-1)} (\ln (t))^{-4q/(p-1)},
\end{equation*}%
for $t$ sufficiently large and for all non vanishing smooth initial data.
\end{example}

\begin{example}\label{ex3.4}[Faster than Polynomials, Less than Exponential]
 Let  $g\left( x\right) =e^{-\left( \ln \left( 1/x\right)
\right) ^{p}}$ on  $(0,r_{0}]$ with  $p>2$. Then the energy of  $\eqref{wave}$ decays as 
\begin{equation*}
E\left( t\right)\leq C_{E\left( 0\right) }e^{-2\left( \ln \left(t\right)\right)^{1/p}}  ,
\end{equation*}%
for $t$ sufficiently large and for all non vanishing initial data in the energy space. 
Moreover the energy of \eqref{wave1D} satisfies the
estimate 
\begin{equation*}
E\left( t\right)  \geqslant  C^{\prime}_{E_1(0)}  e^{-4\left( \ln \left(t\right)\right)^{1/p}}  ,
\end{equation*}%
for $t$ sufficiently large and for all non vanishing smooth initial data.

\end{example}

\begin{rem}
For these above four examples, one can show that $\lim\limits_{x\rightarrow 0^{+}}\Lambda _{H}\left( x\right)<1$. 
Moreover, it is proved in~\cite{AB10}, that the above decay rates are optimal in the finite dimensional case.
\end{rem}

%
%
%

 \section{Weighted nonlinear integral inequalities and decay rates}\label{sec4}

\begin{defn}
We say that a nonnegative function $E$ defined on $[0,+\infty)$ satisfies a weighted nonlinear integral inequality if
there exists a nonnegative function $w$ defined on $[0, \eta)$ with $0<\eta\leq+\infty$  and 
a constant $M>0$ such that $E ([0,+\infty)) \subset [0,\eta)$ and
\begin{equation}\label{NLII}
\int_t^{\infty} w(E(s))E(s) ds\leq M E(t) \quad \forall \ t  \geqslant  0 \,.
\end{equation}
\end{defn}
\begin{defn}
We say that such a weight function $w$ is of "polynomial" type,  if there exists $\alpha>0$ such that
$$
w(s)= s^{\alpha}, \quad \forall \ s  \geqslant  0, \mbox{ with } \eta=+\infty.
$$
\end{defn}
\noindent It is well-known that when $E$ is a nonnegative, nonincreasing absolutely continuous function satisfying
\eqref{NLII} with a polynomial weight function, then $E$ satisfies an optimal decay rate at infinity, proved in \cite[Theorem 9.1] {kom} (see also herein for other references),
that we recall in the next subsection. 

\subsection{Polynomial weights}
\begin{lem}\cite[Theorem 9.1] {kom}  \label{lem-poly}
 Assume that $ E: [0,+\infty) \mapsto [0,+\infty)$ is a non-increasing function and  that there are two
constants $\alpha>0$ and $T>0$ such that
  $$ \int_t^\infty E^{\alpha+1}(s)ds\leq  TE(0)^\alpha E(t), \quad\forall t\in [0,+\infty).$$
Then we have 
  $$E(t)\leq  E(0) \left(\f{T+\alpha t}{T+\alpha T} \right) ^{-\f{1}{\alpha}},\quad\forall t\in [T,+\infty).$$
\end{lem}
\begin{rem}
Note that when applied to decay estimates for dissipative systems, the above Lemma has to be used
with a constant $T$ which blows up as $E(0)=0$, so that the minimal time for which the above decay estimate is valid,
blows up as $E(0)=0$. This result can be reformulated as below, to give an estimate
which is valid for $E(0)>0$ as well as for $E(0)=0$ and for any $t  \geqslant  0$ as explained below.
\end{rem}

\begin{cor}\label{decayPoly}
Assume that $ E: [0,+\infty) \mapsto [0,+\infty)$ is a non-increasing function and  that there are two
constants $\alpha>0$ and $M>0$ such that
 $$ \int_t^\infty E^{\alpha+1}(s)ds\leq  M E(t), \quad\forall t\in [0,+\infty).$$
Then we have 
$$
E(t)\leq  E(0) \min\Big(\Big(\frac{M(\alpha +1)}{M+\alpha E(0)^{\alpha} t}\Big)^{1/\alpha}, 1\Big) ,\quad\forall \  t  \geqslant  0.
$$
\noindent In particular for $E(0)>0$, we deduce that
$$
E(t)\leq E(0) \Big(\frac{M(\alpha +1)}{M+\alpha E(0)^{\alpha} t}\Big)^{1/\alpha},
 \quad \forall \ t  \geqslant  M E(0)^{-1/\alpha} \,.
$$
\end{cor}

\begin{lem}\cite[Theorem 8.1] {kom}  \lab{lem-expo}
Assume that $ E: [0,+\infty) \mapsto [0,+\infty)$ is a non-increasing function and that  there is a constant  $T>0$ such that
   $$\int_t^\infty E(s)ds\leq  T E(t), \quad\forall t\in  [0,+\infty).$$
Then we have
   $$E(t)\leq  E(0)e^{1-t/T},\quad\forall  t\in [T,+\infty).$$
\end{lem}

\subsection{General weights}

For general weight functions, semi-explicit optimal decay rates have been derived for the first time in \cite{AB05}, and later on a simplified
form of the rates in \cite{AB10}.

Let $\eta >0$ and $M>0$ be fixed  and $w$
be a strictly increasing function from $[0,\eta)$ onto $[0,+\infty)$.
For any $r \in (0, \eta)$, we define a function $K_r$ from $(0,r]$ on $[0,+\infty)$
by:
\begin{equation}\label{K}
K_r(\tau)=\int_{\tau}^{r} \frac{dy}{y w(y)} \,,
\end{equation}
\noindent and a function
$\psi_r$ which is a strictly increasing onto function defined
from $[\frac{1}{w(r)},+\infty)$ on $[\frac{1}{w(r)}, +\infty)$ by:
\begin{equation}\label{psir}
\displaystyle{\psi_r(z)=z+K_r(w^{-1}(\frac{1}{z})) \geq z  \ , \quad \forall \ z \geq \frac{1}{w(r)}\,.}
\end{equation}

\noindent We can now formulate our weighted integral inequality:

\begin{thm}\label{thmint}\cite[Theorem 2.1] {AB05}
We assume that $E$ is a nonincreasing, absolutely continuous
 function from $[0,+\infty)$ on $[0,\eta)$, satisfying
the inequality
\begin{equation}\label{ineqint}
\int_{S}^T w(E(t)) E(t)  \, dt \leq M E(S) \,, \quad \forall \, 0 \leq S \leq T \,.
\end{equation}
\noindent Then $E$ satisfies the following estimate:
\begin{equation}\label{decayenergyamo}
\displaystyle{E(t) \leq w^{-1}\Big(\frac{1}{\psi_r^{-1}(\frac{t}{M})}\Big) \ , \quad \forall \ t\geq
\frac{M}{w(r)}}\,,
\end{equation}
where $r>0$ is such that
$$
\frac{1}{M}\int_0^{+\infty} E(\tau)w(E(\tau))\, d\tau \leq r
\leq \eta \,.
$$

\noindent In particular, we have $\displaystyle{\lim_{t \rightarrow +\infty}E(t)=0}$ with  the decay rate
given by  (\ref{decayenergyamo}).
\end{thm}

\begin{thm}\label{thmintc} \cite[Theorem 2.3]{AB10}
Let $H$ be a strictly convex function on $[0,r_0^2]$ such the $H(0)=H'(0)=0$. We define $L$ and $\Lambda_H$ as above.
Let $E$ be a given nonincreasing, absolutely continuous
 function from $[0,+\infty)$ on $[0,+\infty)$ with $E(0)>0$, $M>0$ and $\beta$ is a given parameter such that
 $$
 \frac{E(0)}{2L(H'(r_0^2))}\leq\beta.
 $$
In addition,  $E$ satisfies the following weighted nonlinear inequality
\begin{equation}\label{ineqint1}
\int_{S}^T L^{-1}(\frac{E(t)}{2\beta}) E(t)  \, dt \leq M E(S) \,, \quad \forall \, 0\leq S \leq T \,.
\end{equation}
\noindent Then $E$ satisfies the following estimate:
\begin{equation}\label{decayenergy}
\displaystyle{E(t)\leq2 \beta L\Big(\frac{1}{\psi_{  0}^{-1}(\frac{t}{M})}\Big) \ , \quad \forall \ t\geq
\frac{M}{H'(r_0^2)}}\,.
\end{equation}
\noindent Furthermore, if \, $\limsup \limits_{x \rightarrow 0^+} \Lambda_H(x)<1$, then $E$ satisfies the following simplified decay rate
\begin{equation}\label{decayenergysimp}
\displaystyle{E(t)\leq2 \beta \Big(H'\Big)^{-1}\Big(\frac{\kappa M}{t}\Big) \ , }
\end{equation}
\noindent for $t$ sufficiently large, and where $\kappa >0$ is a constant independent of $E(0)$.
\end{thm}

\section{Proof of main results} \label{sec5}

In this section, we prove the main results including Theorems \ref{thm-wellposed}-\ref{thm-wave-lower} and the decay rates in Examples \ref{ex3.1}-\ref{ex3.4}.

\subsection{Proof of Theorem \ref{thm-wellposed}}
We define the following unbounded nonlinear operator $\mathcal{A}$ in $\mathcal{H}$ by
$$
\mathcal{A}U=(p, \Delta u -\alpha q - \rho(.,p), q, \Delta v +\alpha p) \,,
$$
with the domain
$$
D(\mathcal{A})=\{U \in \mathcal{H}; \ \mathcal{A}U \in \mathcal{H}\} \,.
$$
It is easy to check that $D(\mathcal{A})=  H^2  \cap ( H_0^1 )\times  H_0^1 )^2$. Moreover, since
$\rho$ is nondecreasing with respect the second variable, we have for all $U, \tilde U \in D(\mathcal{A})$,
$$
\langle \mathcal{A}U - \mathcal{A}\widetilde{U} \,, U - \widetilde{U} \rangle=- \int_{\Omega} (\rho(x,p) - \rho(x, \widetilde{p}))(p
- \widetilde{p})dx\leq0 \,,
$$
Thus $-\mathcal{A}$ is a monotone operator. We now claim that $-\mathcal{A}$ is a maximal operator. We proceed as follows. We denote by $A$ the unbounded operator in $  L^2 $ defined by $A=-\Delta$ and $D(A)=  H^2 \cap  H_0^1 $. Then $I-A$ is invertible as an operator acting from $ H_0^1 $ in $H^{-1}(\Omega)$, so that the operator
$(I-A)^{-1}$ is a well-defined, self-adjoint and if $w \in   L^2 $ then $(I-A)^{-1}w \in   H^2  \cap  H_0^1 $.
Then for any $F=(f,g,h,r) \in \mathcal{H}$, the equation 
$$
(I- \mathcal{A})U=F
$$
with $U=(u,p,v,q) \in D(\mathcal{A})$ is equivalent to
\begin{equation}\label{uu}
\begin{cases}
u - \Delta u + \alpha (I-A)^{-1}(\alpha u) + \rho(., u-f)=G  \,, \\
v=(I-A)^{-1}(H_2 + \alpha u) \,, \\
p= u-f \,, q= v-h \,,
\end{cases}
\end{equation}
where 
\begin{equation}\label{hh}
\begin{cases}
H_1= g+f + \alpha h \in   L^2  \,, \\
H_2= r+h - \alpha f \in   L^2  \,, \\
G= H_1 - \alpha (I-A)^{-1}H_2 \in   L^2 \,.
\end{cases}
\end{equation}
We define for $\theta \in \mathbb{R}$
$$
R(x, \theta)=\int_0^{\theta} \rho(x,s) ds \,.
$$
Let us define the functional $J:  H_0^1  \mapsto \mathbb{R}$ defined by
$$
J(u)= \int_{\Omega} \Big(\frac{1}{2}(u^2 + |\nabla u|^2 + \Big|(I-A)^{-1/2}(\alpha u)\Big|^2) + R(x,u-f) -G u\Big)dx \,.
$$
Note that thanks to our hypotheses, $|\rho(x,s)|\leq C(1+|s|)$ for all $(x,s) \in \Omega \times \mathbb{R}$, so that $J$
is well-defined and continuously differentiable on $ H_0^1  $. Moreover, we have
$$
J^{\prime}(u).\varphi= \int_{\Omega} (u \varphi + \nabla u \cdot \nabla \varphi  + (I-A)^{-1/2}(\alpha u)\,
(I-A)^{-1/2}(\alpha \varphi) + \rho(x,u-f)\varphi -G \varphi)dx\,.
$$
We denote by $||\cdot||$ the $ L^2 $ norm. Since $\rho$ is nondecreasing with respect to the second variable, $J$ is a convex function and we also have
$$
J(u)  \geqslant  \Big( \frac{||u||}{2} -  ||G||\Big)||u|| + \frac{||\nabla u||^2}{2} \,,
$$
so that $J(u) \longrightarrow +\infty$ as $||\nabla u|| \longrightarrow +\infty$. Hence $J$ is coercive. Therefore $J$ attains
a minimum at some point $u \in  H_0^1 $, which satisfies the Euler equation
$$
J^{\prime}(u)=0 \,.
$$
The usual elliptic theory implies that the weak solution $u$ of the variational problem
$$
\begin{cases}
u \in  H_0^1  \,,\\
J^{\prime}(u).\varphi= 0  \,, \quad \forall \  \varphi \in  H_0^1
\end{cases}
$$
is in $  H^2 $. Hence $u \in   H^2  \cap  H_0^1 $. By defining $v$ as in \eqref{uu}, and $p\,,q$
as in \eqref{uu}, it follows that  $U=(u,p,v,q) \in D(\mathcal{A})$ and  $(I- \mathcal{A})U=F$. Hence
$-\mathcal{A}$ is a maximal monotone operator. 
We conclude Theorem \ref{thm-wellposed} using the classical theory of maximal monotone operator (see e.g. \cite{kom} and the references therein).

\subsection{Proof of Theorems \ref{thm-wave} and \ref{thm-wave-gen_decay}}

The proof will be divided in three steps, following those described in \cite{ABCime} (see also \cite{AB05, AB10}).

\begin{description}
\item {\bf Step 1:} we first prove that the energy $E$ satisfies a suitable dominant energy estimate. This is the step in which the geometric assumptions $PMGC$ on both the damping and the coupling regions are used, together with suitable multipliers adapted to the coupled structure of the wave-wave system. The proof is valid without specifying the growth assumptions on the feedback $\rho$.

\item {\bf Step 2:} we then prove that nonnegative and nonincreasing functions $E$ satisfying a suitable dominant energy estimate, satisfies a general weighted nonlinear inequality. 
In the case of polynomially growing feedbacks $\rho$, the proof is easier since the weight function for integral inequalities is known. The general growing case relies on the optimality-convexity method of the first author \cite{AB05}.

\item {\bf Step 3:} we deduce energy decay rates, applying Corollary \ref{decayPoly} for polynomially growing feedbacks, whereas applying Theorem~\ref{thmintc} for general growing feedbacks.

%
%

\end{description}

Let us start with {\small {\bf Step 1.}} 
We use the dominant energy method as developed and explained by the first author in \cite{AB05, ABCime}. This method consists
in estimating time integrals of the nonlinear weighted energy of the system by corresponding dominant weighted energies, here in the frictional case, it means by respectively the nonlinear kinetic energy and the localized linear kinetic energy. Note that this step is valid for feedbacks with polynomial as well as arbitrary growth close to the origin.

\begin{thm} \label{thm-wave-dominant-energy}[Weighted dominant energy method]
We assume that \eqref{HC} holds where $\omega_c$ satisfies the $PMGC$ and that $\rho\in C(\overline{\Omega} \times \mathbb{R})$ is nondecreasing with respect to the second variable, and $\rho(\cdot,0) \equiv 0$ on $\Omega$. Let 
$\omega_d$ be a given subset of $\Omega$ satisfying the $PMGC$.  
 Let $\phi: [0,+\infty) \mapsto [0,+\infty)$ be a non-increasing and absolutely continuous function.
Then, there exist constants 
$\delta_i>0 (i=1,2,3)$ and $\alpha^{\ast}>0$ depending only on $\Omega\,, \omega_d$ but independent of $\phi$
such that 
for any  initial data in the energy space, 
for all $\al^+ \in (0,\alpha^{\ast}]$,  the total energy of the system \eqref{wave}
satisfies the following nonlinear weighted estimate
 \be \lab{decay-dominant-energ} 
  \begin{array}{rl}     \displaystyle \int_S^T \phi(t)E(t)dt 
   \leq  &  \delta_1\phi(S)E(S) +\\
   &\displaystyle \delta_2\int_S^T\phi(t) \int_\Omega \rho^2(x,u')dxdt  
         +\delta_3\int_S^T\phi(t)\int_{\omega_d} |u'|^2dxdt,
     \end{array}
     \ee
 \end{thm}
 
Before proving Theorem \ref{thm-wave-dominant-energy}, we give a Lemma on a weighted energy estimate 
for a non-homogeneous wave equation. The proof of Lemma \ref{lem-wave} will be given in Section \ref{sec6}.

\begin{lem} \label{lem-wave}
Let $\omega$ be a nonempty open subset of $\Omega$ satisfying the PMGC.  Let $\phi: [0,+\infty) \mapsto [0,+\infty)$ be a non-increasing and absolutely continuous function.
Then, there exist constants 
$\eta_i>0 \ (i=1,2,3,4)$ independent on $\phi$ such that 
for all $(u^0,u^1) \in  H_0^1 \times   L^2  $, 
all $f\in  L^2((0,\infty);  L^2(\Omega) )$ and all $0\leq  S\leq  T$, 
the solution $u$ of
 \be\lab{eqn-wave}
 \left\{\begin{array}{l} 
   u''-\Delta u = f \quad\quad   \text{in} \ \Omega\times (0,+\infty),\\
   u= 0\quad  \text{ on } \  \Gamma\times  (0,+\infty),\\
  (u,u')(0)=(u^0,u^1) \quad\quad\ \text{in} \ \Omega 
         \end{array} 
         \right.
    \ee 
satisfies the estimate
 \be\lab{est-lemma}
 \begin{array}{rl}
 &\displaystyle \int_S^T\phi(t)e(t)dt\\
 \leq &  \displaystyle  \eta_1 \phi(S)[e(S)+e(T)]+\eta_2\int_S^T-\phi'(t)e(t)dt 
    \\
 &\displaystyle+\eta_3\int_S^T\phi(t)\int_\Omega  |f|^2 dxdt
   +\eta_4\int_S^T\phi(t)\int_{\omega}  |u'|^2dxdt.
 \end{array}
 \ee
where $e(t):=\f12\int_\Omega(|u'|^2+|\nabla u|^2)dx.$  
 \end{lem}

\noindent
{\bf Proof of Theorem \ref{thm-wave-dominant-energy}.} 

We first consider smooth initial data, then
system \eqref{wave} admits a unique solution $(u,v) \in  C([0,+\infty),(  H^2  \cap  H_0^1 )^2)\cap
W^{1,\infty}([0,+\infty),( H_0^1 )^2)$. 

 Let the weight function $\phi$ be  a non-increasing absolutely continuous function. Let then
 $$e_1(t):=\f12\int_\Omega(|u'|^2+|\nabla u|^2)dx. $$
\noindent We now apply Lemma \ref{lem-wave} to the first equation of \eqref{wave} with $\omega=\omega_d$, $f=-\rho(.,u')-\alpha(x)v'$, with $E$ given by \eqref{ene-wave} 
and $e(t)=e_1(t)$. 
 Using  $e_1 (t)\leq E(t)$ and the property that $E$ is nonincreasing, we obtain for all $0\leq  S\leq  T$ and some constants $\eta_i>0 \ (i=1,2,3,4)$ that
 \be\lab{b1}
 \begin{array}{rl}&  \displaystyle  \int_S^T\phi(t)e_1(t)dt   \\ 
 \leq &\displaystyle  \eta_1\phi(S)[e_1(S)+e_1(T)]+\eta_2\int_S^T(-\phi'(t))e_1(t)dt\\
  &\displaystyle+2\eta_3\int_S^T\phi(t)\int_\Omega\alpha^2(x)|v'|^2dxdt
      + 2\eta_3\int_S^T\phi(t)\int_\Omega   \rho^2(x,u')dxdt   \\
  &\displaystyle+\eta_4\int_S^T\phi(t)\int_{\omega_d}  |u'|^2dx dt  \\ 
  \leq &\displaystyle   C_1\phi(S)E(S)+ 2\eta_3\int_S^T\phi(t)\int_\Omega   \rho^2(x,u')dxdt   \\
  &\displaystyle+\eta_4\int_S^T\phi(t)\int_{\omega_d} |u'|^2dxdt  
         +2\eta_3\int_S^T\phi(t)\int_\Omega\alpha^2(x)|v'|^2dxdt\,.
   \end{array}
   \ee

We now set
$$
e_2(t):=\f12\int_\Omega(|v'|^2+|\nabla v|^2)dx\,,
$$
\noindent and secondly, we apply Lemma \ref{lem-wave} to the second equation of \eqref{wave} with $\omega=\omega_c$, $f=\alpha(x)u'$, with $E$ given by \eqref{ene-wave}, $e(t)=e_2(t)$.  Again, using the inequality $e_2(t)\leq  E(t)$ and the property that $E$ is nonincreasing, we obtain for all $0\leq  S\leq  T$ and some constants $\gamma_i>0 \ (i=1,2,3,4)$, 
 \be\lab{b2}
  \begin{array}{rl} 
&   \displaystyle  \int_S^T\phi(t)e_2(t)dt  \\ 
  \leq & \displaystyle  \gamma_1\phi(S)[e_2(S)+e_2(T)]+\gamma_2 \int_S^T(-\phi'(t))e_2(t)dt  \\
    &\displaystyle+\gamma_3\int_S^T\phi(t)\int_\Omega \alpha^2(x)|u'|^2dxdt
        + \gamma_4\int_S^T\phi(t)  \int_{\omega_c}|v'|^2dt\\[0.3cm]
    \leq&\displaystyle  C_2\phi(S)E(S) +\gamma_4\int_S^T\phi(t)\int_{\omega_c} |v'|^2dxdt 
                       +\gamma_3\int_S^T\phi(t)\int_\Omega \alpha^2(x)|u'|^2dxdt.
  \end{array}
  \ee

Let $\delta>0$ be a real parameter to be chosen later on. Adding \e{b1} to $\delta \cdot $\e{b2}, 
 we obtain that for all $0\leq  S\leq  T$ and all $\delta>0$
 \be\lab{b3}
 \begin{array}{rl} 
 &\displaystyle \int_S^T\phi(t)\big(e_1(t) + \delta e_2(t)\big)dt\\ 
\leq&\displaystyle C(1+\delta)\phi(S)E(S)+ C\int_S^T\phi(t) \int_\Omega  \rho^2(x,u')dxdt\\
&\displaystyle + C\int_S^T\phi(t)\int_{\omega_d}  |u' |^2dx dt
 + C\big( \alpha_+ + \frac{\delta}{\alpha_-}\big) \int_S^T\phi(t)\int_\Omega \alpha |v' |^2dxdt\\
&\displaystyle + C\delta \alpha_+\int_S^T\phi(t)\int_\Omega \alpha|u'|^2dxdt \,,
  \end{array}
  \ee 
where $C$ denotes generic positive constants which may vary from one line to another.
Next, we estimate the term $\int_S^T \phi(t)\int_\Omega \alpha(x)|v'|^2dxdt$ through the coupling relation. 
Obviously, the following identity holds for the solution $(u,v)$ of system \eqref{wave}: 
   $$
  \int_S^T\phi(t)\int_\Omega [ v' (u''-\Delta u +\alpha(x)v' +\rho(x,u'))  
             + u'  (v''-\Delta v-\alpha(x)u')  ]dxdt=0.
   $$ 
After integration by parts, we obtain by Cauchy-Schwartz inequality that for all $\varepsilon_1>0$
 \be\lab{est-v'-wave}
 \begin{array}{rl}
 &\displaystyle \int_S^T\phi(t)\int_\Omega \alpha(x)|v'|^2dxdt   \\ 
 =&\displaystyle   \int_S^T\phi(t)\int_\Omega \alpha(x)|u'|^2dxdt  
           -\left [  \phi(t)\int_\Omega(u'v'+\nabla u\cdot \nabla v)dx  \right]^T_S  \\
 &\displaystyle  +\int_S^T \phi'(t)\int_\Omega (u'v'+\nabla u\cdot \nabla v)dxdt  \\
 &\displaystyle  - \int_S^T\phi(t) \int_\Omega \rho(x,u')v'dxdt   \\
\leq&\displaystyle \int_S^T\phi(t)\int_\Omega \alpha |u'|^2dxdt    +C\phi(S)E(S)    \\
   &\displaystyle + \varepsilon_1 \int_S^T\phi(t)  \int_\Omega|v'|^2 dxdt 
   +\f{C}{\varepsilon_1} \int_S^T\phi(t)  \int_\Omega \rho^2(x,u')dx dt.
   \end{array}
  \ee
Using \eqref{est-v'-wave} in \eqref{b3}, we obtain for all $0<\delta\leq1$ and all $\varepsilon_1>0$ that
 \be\lab{newb3}
 \begin{array}{rl} 
 &\displaystyle \int_S^T\phi(t)\big(e_1(t) + \delta e_2(t)\big)dt\\ 
\leq&\displaystyle C\big(1+\delta + \alpha_+ +\frac{\delta}{\alpha_-}\big)\phi(S)E(S)+ 
 C\Big(1 + \big( \alpha_+ + \frac{\delta}{\alpha_-}\big) \frac{1}{\varepsilon_1}\Big)
 \int_S^T\phi(t) \int_\Omega  \rho^2(x,u')dxdt\\
&\displaystyle + C\int_S^T\phi(t)\int_{\omega_d}  |u' |^2dx dt
 + C_2\big( \alpha_+ + \frac{\delta}{\alpha_-}\big)\varepsilon_1 \int_S^T\phi(t)\int_\Omega  |v' |^2dxdt\\
&\displaystyle + C_1\alpha_+\big( \alpha_+ + \frac{\delta}{\alpha_-}\big) \int_S^T\phi(t)\int_\Omega |u'|^2dxdt \,,
  \end{array}
  \ee 
where $C, C_1$ and $C_2$ are generic positive constants.
Thus, 
 \be\lab{newb4}
 \begin{array}{rl} 
 &\displaystyle 
 \big(1-2C_1\alpha_+\big(\alpha_+ +\frac{\delta}{\alpha_-}\big)\big)\int_S^T\phi(t)e_1(t) dt +
  \big(\delta-2C_2\varepsilon_1\big(\alpha_+ +\frac{\delta}{\alpha_-}\big)\big)\int_S^T\phi(t)e_2(t) dt\\ 
\leq&\displaystyle C\big(1+\delta + \alpha_+ +\frac{\delta}{\alpha_-}\big)\phi(S)E(S)+ 
 C\Big(1 + \big( \alpha_+ + \frac{\delta}{\alpha_-}\big) \frac{1}{\varepsilon_1}\Big)
 \int_S^T\phi(t) \int_\Omega  \rho^2(x,u')dxdt\\
&\displaystyle + C\int_S^T\phi(t)\int_{\omega_d}  |u' |^2dx dt \,.
  \end{array}
  \ee 
Let $\alpha_+$ be small so that 
$$
0<\alpha_+\leq\alpha^{\ast}:=\sqrt{\frac{1}{4 C_1}} \,.
$$
We then fix $\delta>0$ so that
$$
\big(1-2C_1\alpha_+\big(\alpha_+ +\frac{\delta}{\alpha_-}\big)\big)  \geqslant  \frac{1}{2} \,
$$
and choose next  $\varepsilon_1>0$ so that
$$
\delta-2C_2\varepsilon_1\big(\alpha_+ +\frac{\delta}{\alpha_-}\big)\big)  \geqslant  \frac{\delta}{2} \,.
$$
With these successive choices of $\alpha_+, \delta$ and $\varepsilon_1$, we deduce that
 \be\lab{newb5}
 \begin{array}{rl} 
 \displaystyle 
\int_S^T\phi(t)E(t) dt\leq 
& C\big(1+\delta + \alpha_+ +\frac{\delta}{\alpha_-}\big)\phi(S)E(S)\\
&\displaystyle + 
 C\Big(1 + \big( \alpha_+ + \frac{\delta}{\alpha_-}\big) \frac{1}{\varepsilon_1}\Big)
 \int_S^T\phi(t) \int_\Omega \rho^2(x,u')dxdt\\
&\displaystyle + C\int_S^T\phi(t)\int_{\omega_d}  |u' |^2dx dt \,.
  \end{array}
  \ee 
For initial data in the energy space, we conclude by a standard argument using density of $  (H^2  \cap  H_0^1) \times  H_0^1 $ in $ H_0^1 \times   L^2 $, together with the dissipativity of the underlying nonlinear semigroup. This ends the proof of   Theorem \ref{thm-wave-dominant-energy}.
\qed

\begin{rem}
The constant $C_1$, and therefore the constant $\alpha^{\ast}$ depends only on $\Omega\,, \omega_d$.
\end{rem}

\noindent {\bf Proof of Theorem \ref{thm-wave}.}[Case of polynomially growing feedbacks]

{\small {\bf Step 2.}}
Assume that $(HF_p)$ holds.

We first consider smooth initial data, then
system \eqref{wave} admits a unique solution $(u,v) \in  C([0,+\infty),(  H^2  \cap  H_0^1 )^2)\cap
W^{1,\infty}([0,+\infty),( H_0^1 )^2)$. 
Moreover, the total energy $E$ defined by \eqref{ene-wave} is absolutely continuous 
and is non-increasing due to the monotonicity of  $\rho$, that is
 \be \label{decay}
 E'(t)=-\int_\Omega \rho(x,u')u'dx\leq  0.
 \ee 

Let $t  \geqslant  0$ be fixed and $\omega_d^{0,t}:=\{x\in \omega_d: \  |u'(t,x)| \geqslant  1\}$ and 
$\omega_d^{1,t} :=\{x\in \omega_d:  \  |u'(t,x)|\leq  1\}$.   In short, we just write $\omega_d^0, \omega_d^1$ in the sequel.
Then it follows from  \eqref{HF} and \eqref{decay}  that
  \be\lab{c3}
    \int^T_S\phi(t)\int_{\omega_d^0}|u'|^2dxdt 
       \leq-C\int^T_S\phi(t)E'(t)dt
       \leq C\phi(S)E(S)
      \ee 
Similarly  we obtain from \eqref{HF}, \eqref{decay} and Young inequality that for every $\varepsilon_2>0$,
 \be \lab{c4} 
 \begin{split}
  \int_S^T \phi(t)  \int_{\omega_d^1}|u'|^2dx dt
& \leq  C  \int_S^T \phi(t)\int_\Omega|\rho(x,u')u'|^{\f{2}{p+1}}dx dt   \\ 
   &   \leq C \int_S^T \phi(t) \left( \int_\Omega |\rho(x,u') u'|dx \right)^{\f{2}{p+1}}dt \\
     &  \leq C \int_S^T \phi(t)(-E'(t))^{\f{2}{p+1}}dt \\ 
    &\leq\int^T_S [\varepsilon_2 (\phi(t))^{\frac{p+1}{p-1}} -C(\varepsilon_2)E'(t)]dt  \\
   &\leq   \varepsilon_2\int^T_S     (\phi(t))^{\frac{p+1}{p-1}}dt+C(\varepsilon_2)E(S).
       \end{split}
       \ee
where $C(\varepsilon_2)>0 $ stands for a constant depending on $\varepsilon_2$ (going to $+\infty$
as $\varepsilon_2$ goes to zero).       
Summing \e{c3} and \e{c4} gives 
  \be \lab{c5}
  \int_S^T\phi(t)\int_{\omega_d}|u'|^2dxdt \leq
      \max\Big(C(\varepsilon_2),C\phi(S)\Big)E(S)+\varepsilon_2\int_S^T (\phi(t))^{\frac{p+1}{p-1}}dt,
       \ee

Similarly, let  $\Omega^0:=\{x\in \Omega:\ |u'| \geqslant  1\}$ and 
$\Omega^1:=\{x\in \Omega:\ |u'|\leq  1\}$.  Note that these subsets depend as above on $t$.
We get from \eqref{HF}, \eqref{decay} and Young inequality  that for every $\varepsilon_3>0$
\be \lab{c6} 
   \begin{split}
    \int_S^T \phi(t)\int_\Omega \rho^2(x,u')dxdt  
&  =    \int_S^T \phi(t) \left( \int_{\Omega^0} \rho^2(x,u')dx
  +\int_{\Omega^1} \rho^2(x,u')dx\right)dt   \\
  &\leq C \phi(S)E(S) + \int_S^T \phi(t)\int_\Omega  |u' \rho(x,u')|^{\frac{2}{p+1}}dxdt \\
   &\leq   \max\Big (C(\varepsilon_3),C\phi(S)\Big)E(S)+\varepsilon_3\int_S^T (\phi(t))^{\frac{p+1}{p-1}}dt.
       \end{split}
\ee

We now choose the weight function $\phi$ as follows
\begin{equation}\label{phipoly}
\phi(t)=E^{\frac{p-1}{2}}(t) \,,\quad \forall \ t  \geqslant  0 \,.
\end{equation}
\noindent Combining  \e{decay-dominant-energ} , \e{c5}, \e{c6}, together with this choice for $\phi$,  and letting $\varepsilon_2$, $\varepsilon_3$ small enough,  
we obtain that for  all $0\leq  S\leq  T$, 
$$
\int_S^T E^{\frac{p+1}{2}}(t) dt\leq\max\Big(C_1, C_2 E^\f{p-1}{2}(0)\Big)E(S)\,,
$$
\noindent where $C_1,C_2$ are positive constants independent of $E(0)$.
  
 We finish the proof by {\small {\bf Step 3}}, which says that nonnegative, nonincreasing functions $E$ satisfying a polynomial nonlinear integral inequality, then satisfy a polynomial decay rate. Applying Corollary \ref{decayPoly} with $\alpha=\frac{p-1}{2}>0$, 
 and $M=C(1+E^{\frac{p-1}{2}}(0))$, we end the proof of Theorem \ref{thm-wave} and obtain when $E(0)>0$
\be\lab{bu2} E(t)\leq C_{E(0)} t^{-2/(p-1)}, \quad \forall \ t  \geqslant  C_1E^{\frac{-(p-1)}{2}}(0) + C_2\,,
\ee
\noindent where
$$
C_{E(0)}= \Big(\max\Big(C_1, C_2 E^\f{p-1}{2}(0)\Big) (1 + \frac{1}{\alpha})\Big)^{1/\alpha} \,.    \qed
$$
 

\noindent {\bf Proof of Theorem \ref{thm-wave-gen_decay}.} [Case of general growing feedbacks]

In {\small {\bf Step 2} },   We shall prove the following theorem

\begin{thm} \label{thm-wave-NLW-integral-ineq}[Optimal-weight convexity method]
Assume the hypotheses of Theorem~\ref{thm-wave-dominant-energy}. 
Assume furthermore that $(HF_g)$ holds where $g$ is such that the function $H$ defined in \eqref{defH} is strictly convex on $[0,r_0^2]$, and $g'(0)=0$. We define $L$ by \eqref{defL}. Let the initial data be in the energy space and be non vanishing,
 $(u,v)$ be the solution of \eqref{wave} and $E$ be its energy. Then $E$ satisfies the following nonlinear weighted
integral inequality
\begin{equation}\label{NLIE}
\int_S^T L^{-1}\Big(\frac{E(t)}{2\beta}\Big)E(t) dt\leq M E(S) , \quad \forall \ 0\leq S\leq T \,,
\end{equation}
\noindent where $\beta$ and $M$ are respectively given by 
$$
\beta= \max\Big(C_2, \frac{E(0)}{2L(H'(r_0^2))}\Big) \,,
$$
\noindent and
$$
M=2C_1(1 + H'(r_0^2)) \,,
$$
\noindent with $C_1>0,C_2>0$ depending on $\delta_i\ ( i=1,2,3)$, $\Omega$, $\omega_d$ but independent  of $E(0)$.
 \end{thm}
\begin{rem}
This method is called the optimal-weight convexity method according to the property that the weight function $\phi$ is chosen
in an optimal way by setting
$$
\phi(.)=L^{-1}\Big(\frac{E(.)}{2\beta}\Big)
$$
\noindent thanks to suitable convexity arguments relying both on Jensen and Young's inequalities for an appropriate convex function.
\end{rem}

\begin{proof}
We consider as before smooth initial data, then
the solution $(u,v)$ of \eqref{wave} is in $C([0,+\infty),(  H^2  \cap  H_0^1 )^2)\cap
W^{1,\infty}([0,+\infty),( H_0^1 )^2)$. 
Moreover, the total energy $E$ satisfies the dissipation relation \eqref{decay}. Thanks to Theorem \ref{thm-wave-dominant-energy}, we know that
$E$ satisfies the weighted dominant energy estimate \eqref{decay-dominant-energ}. We shall now use the optimal-weight convexity method of the first author \cite{AB05} to build an optimal weight function $\phi$ to prove that the two terms
$$
 \int_S^T\phi(t) \int_\Omega \rho^2(x,u')dxdt  \,  \text{ and }  \int_S^T\phi(t)\int_{\omega_d} |u'|^2dxdt 
$$
in \eqref{decay-dominant-energ}, are bounded above by the term
$$
C E(S)(1+ \phi(S)) + C   \int_S^T E(t)\phi(t) dt \quad \forall \ 0\leq  S\leq T \,.
$$

We proceed as in \cite{AB05}. Choose a parameter $\varepsilon_0$ sufficiently small, e.g. $\varepsilon_0=min(1,g(r_0))$.

For fixed $t  \geqslant  0$, we define the subset $\Omega_1^t=\{
x \in \Omega \,, |u'(t,x)|\leq\varepsilon_0\}$. Now thanks to \eqref{HFg}, we know that \eqref{gwrho} holds. 
Hence, since $g$ is increasing, we have
\begin{equation}\label{int1}
g\big(\frac{|\rho(x,u'(t,x))|}{K}\big)\leq|u'(t,x)| \,, \mbox{ for a.e } x \in \Omega_1^t \,,
\end{equation}
where $K=c_2 ||a||_{\infty}$, with $||\cdot||_{\infty}$ standing for the $L^{\infty}$ norm.
Now, we can note that parameter $\varepsilon_0$ has been chosen to guarantee the following two properties
\begin{equation}\label{dconv1}
\frac{1}{|\Omega_1^t|} \int_{\Omega_1^t} \frac{|\rho(x,u'(t,x))|^2}{K^2} dx \in [0,r_0^2] \,,
\end{equation}
and
\begin{equation}\label{dconv2}
 \frac{1}{|\Omega_1^t|K} \int_{\Omega_1^t} \rho(x,u'(t,x)) u'(t,x) dx \in [0,H(r_0^2)] \,,
\end{equation}
hold.
Since $H$ has been assumed to be convex on $[0,r_0^2]$ and thanks to \eqref{dconv1}, the Jensen's inequality, and \eqref{int1}, we obtain
\begin{multline*}
H\Big(\frac{1}{|\Omega_1^t|} \int_{\Omega_1^t} \frac{|\rho(x,u'(t,x))|^2}{K^2} dx\Big) \leq 
\frac{1}{|\Omega_1^t|} \int_{\Omega_1^t} H\Big(\frac{|\rho(x,u'(t,x))|^2}{K^2}\Big) dx =\\
\frac{1}{|\Omega_1^t|K} \int_{\Omega_1^t} |\rho(u'(t,x))| \,g\big(\frac{|\rho(x,u'(t,x))|}{K}\big)dx\leq
\frac{1}{|\Omega_1^t|K} \int_{\Omega_1^t} u'(t,x)\rho(x,u'(t,x))dx \,.
\end{multline*}
By \eqref{dconv2}, we deduce that
$$
\int_S^T \phi(t) \int_{\Omega_1^t}|\rho(x,u'(t,x))|^2dx dt\leq\int_S^T K^2 |\Omega_1^t| \, \phi(t)
H^{-1} \Big(\frac{1}{|\Omega_1^t|K} \int_{\Omega_1^t} u'(t,x)\rho(x,u'(t,x))dx\Big) dt \,,
$$
and using further Young's inequality, the dissipation relation \eqref{decay}, 
we obtain
$$
\int_S^T \phi(t) \int_{\Omega_1^t}|\rho(x,u'(t,x))|^2dx dt\leq
K^2 |\Omega| \int_S^T  \, \widehat{H}^{\ast}(\phi(t))dt +  K E(S), \quad \forall \ 0 \leq S \leq T \,.
$$
On the other hand, we prove easily as in the proof of Theorem \ref{thm-wave} that
$$
\int_S^T \phi(t) \int_{\Omega \backslash \Omega_1^t}|\rho(x,u'(t,x))|^2dx dt\leq
K E(S) \phi(S), \quad \forall \ 0\leq S\leq T \,.
$$
Adding these two inequalities, we obtain
\begin{equation}\label{nlikn1}
\int_S^T \phi(t) \int_{\Omega}|\rho(x,u'(t,x))|^2dx dt\leq K^2 |\Omega| \int_S^T  \, \widehat{H}^{\ast}(\phi(t))dt + 
KE(S)(1+ \phi(S)) , \quad \forall \ 0\leq S\leq T \,.
\end{equation}
We now turn to the estimate of the localized weighted linear kinetic energy. Thanks to \eqref{HFg}, we know that \eqref{gwrho} holds.
Choose a parameter $\varepsilon_1$ sufficiently small, e.g. 
$\varepsilon_{1}=\min \{r_{0},g(r_{1})\}$ where $r_{1}$ is defined by
$$
r_{1}^2=H^{-1}\left(\frac{k}{K}H(r_{0}^2)\right),
$$
where $k=c_1 a_-$.
For fixed $t  \geqslant  0$, we define the subset $\omega_d^t=\{
x \in \omega_d \,, |u'(t,x)|\leq\varepsilon_1\}$. Thanks to \eqref{gwrho} and since $\eqref{HFg}$ holds, we have
\begin{equation}\label{int2}
g(|u'(t,x)|)\leq\frac{|\rho(x,u'(t,x))|}{k}\,, \mbox{ for a.e } x \in \omega_d^t \,.
\end{equation}
Now, we can note that parameter $\varepsilon_1$ has been chosen to guarantee the following two properties
\begin{equation}\label{dconv3}
\frac{1}{|\omega_d^t|} \int_{\omega_d^t} |u'(t,x)|^2dx \in [0,r_0^2] \,,
\end{equation}
and
\begin{equation}\label{dconv4}
 \frac{1}{|\omega_d^t|k} \int_{\omega_d^t} \rho(x,u'(t,x)) u'(t,x) dx \in [0,H(r_0^2)] \,,
\end{equation}
hold.
Since $H$ has been assumed to be convex on $[0,r_0^2]$ and thanks to \eqref{dconv3}, the Jensen's inequality, and \eqref{int2}, we obtain
\begin{align*}
H\Big(\frac{1}{|\omega_d^t|} \int_{\omega_d^t} |u'(t,x)|^2 dx\Big)
& \leq 
\frac{1}{|\omega_d^t|} \int_{\omega_d^t} H (|u'(t,x)|^2)dx \\
&=\frac{1}{|\omega_d^t|} \int_{\omega_d^t} |u'(t,x)|g(|u'(t,x)|)dx  \\
&\leq 
\frac{1}{|\omega_d^t|k} \int_{\omega_d^t} u'(t,x) \rho(x,u'(t,x))dx \,.
\end{align*}
Thanks to \eqref{dconv4}, we deduce that
$$
\int_S^T \phi(t) \int_{\omega_d^t}|u'(t,x)|^2dx dt\leq\int_S^T |\omega_d^t| \, \phi(t)
H^{-1} \Big(\frac{1}{|\omega_d^t|k} \int_{\Omega_1^t} u'(t,x)\rho(x,u'(t,x))dx\Big) dt \,,
$$
and using further Young's inequality, the dissipation relation \eqref{decay},
we obtain
$$
\int_S^T \phi(t) \int_{\omega_d^t}|u'(t,x)|^2dx dt\leq
|\omega_d| \int_S^T  \, \widehat{H}^{\ast}(\phi(t))dt +  \frac{1}{k}E(S), \quad \forall \ 0\leq S\leq T \,.
$$
On the other hand, we prove easily as in the proof of Theorem \ref{thm-wave} that
$$
\int_S^T \phi(t) \int_{\omega_d \backslash \omega_d^t}|u'(t,x)|^2dx dt\leq
\frac{1}{k}E(S) \phi(S), \quad \forall \ 0\leq S\leq T \,.
$$
Adding these two inequalities, we obtain
\begin{equation}\label{nlikn2}
\int_S^T \phi(t) \int_{\omega_d}|u'(t,x)|^2dx dt\leq |\omega_d| \int_S^T  \, \widehat{H}^{\ast}(\phi(t))dt + 
\frac{1}{k}E(S)(1+ \phi(S)) , \quad \forall \ 0\leq S\leq T \,.
\end{equation}
Using \eqref{nlikn1} and \eqref{nlikn2} in the weighted dominant energy estimate \eqref{decay-dominant-energ}, we obtain
\begin{equation}\label{optweig}
\int_S^T \phi(t) E(t) dt\leq C_1 E(S)(1+\phi(S)) +C_2 \int_S^T  \, \widehat{H}^{\ast}(\phi(t))dt , \quad \forall \ 0\leq S\leq T \,,
\end{equation}
where the constants $C_1, C_2$ depend only on the $\delta_i$ for $i=1,2,3$ and on $|\Omega|$ and $|\omega_d|$ in an explicit way. In particular, they do not depend on  $\phi$.

Let
\begin{equation}\label{defbeta}
\beta= \max\Big(C_2, \frac{E(0)}{2L(H'(r_0^2))}\Big) \,.
\end{equation}
where $L$ is defined in \eqref{defL}. 
Since $E$ is a nonincreasing function, and thanks to \eqref{LH'}, we have
$$
\frac{E(t)}{2\beta}\leq\frac{E(0)}{2\beta}\leq L(H'(r_0^2)) <r_0^2 \,.
$$
Hence, since $L^{-1}$ is defined from $[0,r_0^2)$ onto $[0,+\infty)$, we can define $\phi$ by
\begin{equation}\label{phigene}
\phi(t)=L^{-1}\Big(\frac{E(t)}{2\beta}\Big) \,,\quad \forall \ t  \geqslant  0 \,.
\end{equation}
By definition of $L$, $\phi$ is a nonnegative, non increasing and absolutely continuous function on $[0,+\infty)$.
We first note that
\begin{equation}\label{Ephi}
\phi(S)\leq H'(r_0^2), \quad \forall \ S  \geqslant  0 \,.
\end{equation}
Then, thanks to our "optimal" choice of the weight function $\phi$ and to the definition of $L$, we have
$$
L(\phi(t))=\frac{E(t)}{2\beta}=\frac{\widehat{H}^{\ast}(\phi(t))}{\phi(t)}, \quad \forall \ t  \geqslant  0 \,.
$$
This implies
$$
C_2 \widehat{H}^{\ast}(\phi(t))\leq\beta \widehat{H}^{\ast}(\phi(t))= \frac{1}{2} \phi(t) E(t), \quad \forall \ t  \geqslant  0 \,.
$$
\noindent Combining this estimate together with \eqref{Ephi} in \eqref{optweig}, we obtain
\begin{equation}\label{intE}
\int_S^T L^{-1}\Big(\frac{E(t)}{2\beta}\Big)E(t) dt\leq M E(S) , \quad \forall \ 0\leq S\leq T \,,
\end{equation}
where
\begin{equation}\label{M}
M=2C_1(1 + H'(r_0^2)) \,.
\end{equation}
\end{proof}
 %
 

 We finish the proof of Theorem \ref{thm-wave-gen_decay} by {\small {\bf Step 3}}. 
Thanks to {\bf Step 2}, $E$ satisfies
the weighted nonlinear integral inequality \eqref{intE}, where $\beta$ is defined by \eqref{defbeta}, $M$ is defined
by \eqref{M}. Hence applying Theorem~\ref{thmintc} with this $\beta$ and $M$ we deduce that $E$ satisfies
the decay rate \eqref{decayenergy} in the general case, and the simplified decay rate \eqref{decayenergysimp} if
$\limsup\limits_{x \rightarrow 0^+} \Lambda_H(x)<1$. This  concludes the proof of Theorem~\ref{thm-wave-gen_decay}. 
\qed

\subsection{Proof of Theorem \ref{thm-wave-lower}}
Thanks to our hypotheses, the simplified upper energy estimate \eqref{decayenergysimp} of Theorem~\ref{thmintc} holds, so that
$E(t)$ converges to $0$ as $t$ goes to infinity. Hence, there exists $T_0  \geqslant  0$ such that
\begin{equation}\label{strong0}
E(t)\leq\Big(\frac{r_0^2}{\gamma_s}\Big)^2 \,, \quad  \forall t\geq T_0, 
\end{equation}
where $\gamma_{s}=4 \sqrt{E_1(0)}$. 
Hence 
\begin{equation}\label{intconv}
\gamma_s  \sqrt{E(t)} \in [0,r_0^2] \ \mbox{ for all }t  \geqslant  T_0 \,.
\end{equation}

On the other hand, thanks to the regularity of $u$ (see \cite{AB11} for details), we have 
$$
|| u' (t,.)||_{L^{\infty}(\Omega)}\leq\gamma_s \sqrt{E(t)}
$$
Using this inequality in the dissipation relation
$$
-E^{\prime}(t)= \int_{\Omega}a(x)   u' g(u') dx \,, \quad t  \geqslant  0 \,,
$$
together with \eqref{intconv}, we deduce that for all $t  \geqslant  T_0$
$$
-E^{\prime}(t)\leq||a||_{L^{\infty}(\Omega)}\frac{2}{\gamma_{s}} \sqrt{E(t)}H(\gamma_s  \sqrt{E(t)}) \,.
$$
Therefore, we have
\begin{equation}\label{infK}
E(t) \geq \Big( \frac{1}{\gamma_s} K^{-1}(||a||_{L^{\infty}(\Omega)}(t-T_0))
\Big)^2 \,, \quad \forall \ t  \geqslant  T_0 \,,
\end{equation}
where $K^{-1}$ denotes the inverse function of $K$ defined by
\begin{equation}\label{defK}
K(\tau)=\int_{\tau}^{z_0} \frac{dy}{H(y)} \,, \quad \forall \tau \in (0, \sqrt{E(T_0)})\,,
\end{equation}
where
\begin{equation}\label{defz0}
z_0=\gamma_s \sqrt{E(T_0)}\,.
\end{equation}
We denote by $z$ the solution of the following ordinary differential equation
\begin{equation}\label{eqz}
z^{\prime}(t) + ||a||_{L^{\infty}(\Omega)}H(z(t))=0 \,, z(0)=z_0 \,.
\end{equation}
Then we have the relation
\begin{equation}\label{Kz}
z(t-T_0) =K^{-1}(||a||_{L^{\infty}(\Omega)}(t-T_0)) \,, \quad \forall \ t  \geqslant  T_0 \,.
\end{equation}
We now use the following comparison Lemma, that we recall for the sake of completeness.
\begin{lem}\label{compar}  \cite[Lemma 2.4]{AB10} 
Let $H$ be a given strictly convex ${\cal C}^1$ function from $[0,r_0^2]$ to $\mathbb{R}$
such that $H(0)=H^{\prime}(0)=0$, where $r_0>0$ is sufficiently small and define $\Lambda_H$ on $(0,r_0^2]$ by
\eqref{defLambdaH}.

Let $z$ be the solution of the 
ordinary differential equation:
\begin{equation}\label{eq-v-comp}
z^{\prime}(t) + \kappa \,H(z(t))=0 \,, \, 
z(0)=z_0 
\quad t \ge 0 \,,
\end{equation}
\noindent where $z_0>0$ and $\kappa>0$ are given. Then $z(t)$
is defined for every $t \ge 0$ and decays to
$0$ at infinity.
Moreover assume that $(HF_l)$ holds.
Then there exists $T_1>0$ such that
for all $R>0$ there exists a constant $C>0$
such that

\begin{equation}\label{comparaison}
z(t) \geq C (H^{\prime})^{-1}\Big(
\frac{R}{t}\Big) ,\quad
\ \forall \ t \ge T_1\,,
\end{equation}
\noindent where $T_1$ is a positive constant.
\end{lem}
We apply this Lemma to the solution $z$ of \eqref{eqz} with $R=1$ and $\kappa=||a||_{L^{\infty}(\Omega)}$. Thus, there exist two constants $T_1>0$
and $C_s>0$ such that
\begin{equation}\label{comparaison}
z(t)  \geq C_s  (H^{\prime})^{-1}\Big(
\frac{1}{t}\Big) \,,\quad
\ \forall \ t \ge T_1\,,
\end{equation}
Combining \eqref{comparaison} together with \eqref{infK} and \eqref{Kz}, we obtain the lower estimate \eqref{carac-low-wave}.
\qed
\begin{rem}
The constant $C$ of the above Lemma~\ref{compar} depends explicitly on $\kappa$, $R$ (and in addition of $\mu$ if the second
alternative of $(HF_l)$ holds). This dependence is
given in the proof of Lemma~2.4 in \cite{AB10}.
Moreover, one may assume that $r_0=\infty$ in Lemma \ref{compar}. In this case
the interval $[0,r_0^2]$ becomes $[0,+\infty)$.
\end{rem}
\subsection{Proof of the decay rates given in Examples \ref{ex3.1}-\ref{ex3.4}}

{\bf Example \ref{ex3.1}} We have $H(x)=x^{(p+1)/2}$ for $ x \in [0,r_0^2]$. Thus $H'(x)= \frac{p+1}{2}x^{(p-1)/2}$ and
$H$ is strictly convex on a right neighborhood of $0$. Moreover, $\Lambda_H(x)=\frac{2}{p+1} <1$
for all $x \in [0,r_0^2]$. We easily conclude applying \eqref{decayenergysimp} of Theorem~\ref{thmintc} and Theorem \ref{thm-wave-lower} for the
lower estimate in the one-dimensional case.

{\bf Example \ref{ex3.2}} We have $H(x)=\sqrt{x}e^{\frac{-1}{x}}$ for $ x \in [0,r_0^2]$. Thus $H'(x)= \frac{e^{-1/x}}{\sqrt{x}}(\frac{1}{2} + \frac{1}{x})$, 
and
$H$ is strictly convex on a right neighborhood of $0$.
Moreover, we have $\Lambda_H(x)= \frac{1}{(\frac{1}{2} + \frac{1}{x})}$ for all $x>0$ sufficiently close to $0$, so that $\lim\limits_{x \rightarrow 0+}\Lambda_H(x)=0$. We apply \eqref{decayenergysimp} of Theorem~\ref{thmintc}. So we set $x(t)=(H')^{-1}\Big(\frac{\kappa M}{t}\Big)$. Then one can prove that $x(t)$ is equivalent to $\frac{1}{\ln(t) - \ln(\kappa M)}$ as $t$ goes to $+\infty$. We therefore obtain the desired upper bound, using this equivalence. One can show that the second alternative of $(HF_l)$ holds for any $\mu >1$ (see subsection 7.10 in \cite{AB10}). Thus, we obtain in the same way by Theorem \ref{thm-wave-lower} the lower estimates in the one-dimensional case.

{\bf Example \ref{ex3.3}} We have $H(x)=x^{(p+1)/2} (\ln(1/\sqrt{x}))^{q}$ for $ x \in [0,r_0^2]$. Thus $H'(x)= \frac{1}{2}x^{(p-1)/2} (\ln(1/\sqrt{x}))^{q}
\Big( p+1 - q (\ln(1/\sqrt{x}))^{-1}\Big)$ and
$H$ is strictly convex on a right neighborhood of $0$. Moreover,
$\Lambda_H(x)= \frac{2}{p+1 - q (\ln(1/\sqrt{x}))^{-1}}$ for all $x>0$ sufficiently close to $0$, so that $\lim\limits_{x \rightarrow 0+}\Lambda_H(x)=\frac{2}{p+1}$. We apply \eqref{decayenergysimp} of Theorem~\ref{thmintc}. So we set $x(t)=(H')^{-1}\Big(\frac{\kappa M}{t}\Big)$ and
$y(t)=\Big(\frac{2 \kappa M}{t}\Big)^{2/(p-1)}$. Then one can prove that $\Big(\frac{x(t)}{y(t)}\Big)^{(p-1)/2}(\ln(1/\sqrt{x}))^{q}$ is equivalent to $\frac{1}{p+1}$ as $t$ goes to $+\infty$. On the other hand, computing $\ln(x(t))$ and $\ln(y(t))$, we find that $\ln(x(t))$ is equivalent to $\ln(y(t))$ as $t$ goes to $+\infty$. Using this relation in the previous one, we find that $x(t)$ is equivalent to $D t^{-2/(p-1)} \Big( \ln(t)\Big)^{-2q/(p-1)}$, where $D$ is an explicit positive constant which depends on $\kappa, M, p$ and $q$.
We therefore obtain the desired upper estimate, using this equivalence. We obtain by Theorem \ref{thm-wave-lower} the lower estimates in the one-dimensional case.

{\bf Example \ref{ex3.4}} We have $H(x)=\sqrt{x}e^{-\left( \ln \left( 1/\sqrt{x}\right)
\right) ^{p}}$ for $ x \in [0,r_0^2]$. Thus, we have
$
H'(x)= \frac{1}{2\sqrt{x}}e^{-\left( \ln \left( 1/\sqrt{x}\right)~
\right) ^{p}}\Big(1+p \Big(\ln \Big(\frac{1}{\sqrt{x}}\Big)\Big)^{p-1}\Big)
$, and
$H$ is strictly convex on a right neighborhood of $0$.
Moreover, we have $\Lambda_H(x)= \frac{2}{1+p \big(\ln(1/\sqrt{x})\big)^{p-1}}$ for all $x>0$ sufficiently close to $0$, so that $\lim \limits_{x \rightarrow 0+}\Lambda_H(x)=0$. We apply \eqref{decayenergysimp} of Theorem~\ref{thmintc}. So we set $x(t)=(H')^{-1}\Big(\frac{\kappa M}{t}\Big)$ and
$y(t)=e^{-2(\ln (\frac{t}{\kappa M}))^{1/p}}$.
Then one can prove that $\ln(x(t))$ is equivalent to $\ln(y(t))$ as $t$ goes to $+\infty$. We further set $z(t)= \ln(1/\sqrt{x(t)})$ so that $z(t)$ goes to $+\infty$ as $t$ goes to $+\infty$, then we have
$z^p(t)( 1-\theta(t))= \ln (\frac{t}{2\kappa M})$, where $\theta(t)=z^{1-p}(t) +\ln ( 1 +p z^{p-1}(t))z^{-p}(t)$, so that $\theta(t)$ goes to $0$ as $t$ goes to $+\infty$. Hence we have $x(t)=e^{-2(\ln (\frac{t}{\kappa M}))^{1/p}\frac{1}{(1-\theta(t))^{1/p}}}$.  We can check that $\ln (\frac{t}{2\kappa M})^{1/p} (1 - (1- \theta(t))^{1/p})$ goes to $0$ as $t$ goes to $+\infty$. Hence, $x(t)$ is equivalent to $e^{-2(\ln (t))^{1/p}}$ as $t$ goes to $+\infty$.
We therefore obtain the desired upper estimate. One can show that the second alternative of $(HF_l)$ holds for any $\mu >1$ (see subsection 7.10 in \cite{AB10}). Thus, we obtain by Theorem \ref{thm-wave-lower} the lower estimates in the one-dimensional case.

\section{Proof of Lemma \ref{lem-wave}}
\label{sec6}
In this section, we prove Lemma \ref{lem-wave} by the piecewise multiplier method 
which relies on the geometric assumptions PMGC on the subset $\omega \subset \Omega$.
Denoting by $\Omega_j$  and $x_j\ (j=1,\cdots,J)$ the sets and the points given by PMGC, 
 we have $\omega \supset N_{\varepsilon}(\cup_{j=1}^J\gamma_j(x_j)\cup(\Omega\setminus
\cup_{j=1}^J\Omega_j))\cap\Omega.$  Here,
$N_{\varepsilon}(U)=\{x\in \mathbb{R}^n, d(x,U)\leq  \varepsilon\}$ with
$d(\cdot,U)$ the usual euclidian distance to the subset $U$ of
$\mathbb{R}^n$, and $\gamma_j(x_j)=\{x\in \Gamma_j, (x-x_j)\cdot \nu_j>0\},$
where $\nu_j$ denotes the outward unit normal of the boundary
$\Gamma_j=\partial\Omega_j.$

Let $0< \varepsilon_0<\varepsilon_1<\varepsilon_2<\varepsilon$ and
define $ Q_i :=N_{\varepsilon_i}[\cup_{j=1}^J\gamma_j(x_j)\cup(\Omega\setminus \cup_{j=1}^J\Omega_j)] (i=0,1,2)$.  
Since $(\overline{\Omega_j}\backslash  Q_1 )\cap\overline{Q_0}=\emptyset$, 
we introduce a cut-off function $\psi_j\in C_0^\infty(\mathbb{R}^N)$  satisfying
\be \label{psi j}
  0\leq  \psi_j\leq  1,\quad\quad\psi_j=1\quad\text{on}\
\overline{\Omega_j}\backslash Q_1;   \quad\quad\ \psi_j=0 \quad\ \text{on}\ Q_0.
\ee
For $m_j(x)=x-x_j$, we define the $C^1$ vector field on $\Omega$:
\be\lab{h}
h(x)=\left\{\begin{array}{l} \psi_j(x)m_j(x)\quad\quad\text{if}\ x\in  \Omega_j,   j=1,\cdots,J
   \\
0\quad\quad\quad\ \text{if}\ x\in \Omega\backslash \cup_{j=1}^J  \Omega_j\end{array}\right.
\ee 
Using the multiplier $\phi(t)h(x)\cdot \nabla u$ to equation \e{eqn-wave}:
$$\int_S^T \phi(t)\int_{\Omega_j} h(x)\cdot \nabla
u(u''-\Delta u-f)dxdt=0
$$
leads to
\be\lab{5.3}  
  \begin{array}{rl} 
   &\displaystyle \int_S^T \phi(t)  \int_{\Gamma_j} \f{\partial u}{\partial \nu_j}h\cdot \nabla u
       +\f12(h\cdot \nu_j)(|u'|^2-|\nabla u|^2)d\Gamma dt  
   \\  
   = &\displaystyle \left[\phi(t)\int_{\Omega_j}u'h\cdot \nabla u dx\right]_S^T
   -\int_S^T\phi'(t)\int_{\Omega_j}u' h\cdot \nabla u dx dt
   \\
    &\displaystyle+\int_S^T\phi(t)\int_{\Omega_j}\f12 div h(|u'|^2-|\nabla  u|^2)dxdt
    \\ 
    &\displaystyle  +\int_S^T\phi(t)\int_{\Omega_j}\left [\sum_{i,k}\f{\pa u}{\pa x_i}\f{\pa  u}{\pa x_k}\f{\pa h_k}{\pa x_i}
     -h\cdot \nabla  u \, f\right ]dxdt.
     \end{array}
     \ee 
 Thanks to the choice of $\psi_j$, 
 the terms  in the left hand side of \e{5.3} vanish except  on the boundary
  $(\Gamma_j\setminus \gamma_j(x_j))\cap\Gamma$. 
Since $u=0$ on this part of boundary, then $u'=0, \nabla u=\frac{\partial u}{\partial\nu_j}\,\nu_j$. 
  Hence, the left side of \e{5.3} becomes 
 $$\f12\int_S^T \phi(t)\int_{(\Gamma_j\backslash \gamma_j(x_j))\cap\Gamma}
    \Big|\f{\partial u}{\partial \nu_j}\Big|^2 \psi_j m_j\cdot \nu_jd\Gamma dt\leq  0.$$ 
Therefore \e{5.3} implies that 
$$\begin{array}{rl}
 &\displaystyle\int_S^T\phi(t)\int_{\Omega_j}\f12 div h(|u'|^2-|\nabla  u|^2)dxdt
   +\int_S^T\phi(t)\int_{\Omega_j}\sum_{i,k}\f{\pa u}{\pa x_i}\f{\pa  u}{\pa x_k}\f{\pa h_k}{\pa x_i}dxdt 
  \\ 
 \leq&  \displaystyle -\left[\phi(t)\int_{\Omega_j}u'h\cdot \nabla u  dx\right]_S^T
  +\int_S^T\phi'(t)\int_{\Omega_j}u' h\cdot \nabla u dx dt
  \\
  &\displaystyle  +\int_S^T\phi(t)\int_{\Omega_j}h\cdot \nabla  u \,f dxdt\,.\end{array}    
   $$
Since, moreover $h(x)=m_j(x)$ on $\bar{\Omega}_j\backslash Q_1$, we obtain that
$$\begin{array}{rl}
   &\displaystyle \int_S^T\phi(t)\int_{\Omega_j\backslash  Q_1}\f N2 (|u'|^2-|\nabla  u|^2)+|\nabla  u|^2 dxdt
   \\ 
  \leq&\displaystyle  -\left[\phi(t)\int_{\Omega_j}u'h\cdot \nabla u dx\right]_S^T
   + \int_S^T\phi'(t) \int_{\Omega_j}u' h\cdot \nabla u dx dt
   \\
  &\displaystyle
  +\int_S^T\phi(t)\int_{\Omega_j}h\cdot \nabla  u \,fdxdt
  \\
  &\displaystyle -\int_S^T\phi\int_{\Omega_j\cap Q_1}\left [\f12 div h(|u'|^2-|\nabla  u|^2)
  +\sum_{i,k}\f{\pa u}{\pa x_i}\f{\pa  u}{\pa x_k}\f{\pa h_k}{\pa x_i}\right ]dxdt\,.\end{array}
  $$
  Summing the above inequality on $j$ and using the facts  that $\Omega\backslash
  Q_1=\cup_{j=1}^J (\Omega_j)\backslash Q_1$ and $h(x)=0$ on $\Omega\setminus\cup_{j=1}^J \Omega_j$, 
  we obtain
   \be\lab{5.4}
    \begin{array}{rl}
    &\displaystyle \int_S^T\phi(t)\int_{\Omega\backslash Q_1}\f N2 (|u'|^2-|\nabla   u|^2)+|\nabla  u|^2 dxdt
    \\ 
   \leq&\displaystyle  -\left[\phi(t)\int_{ \Omega }u'h\cdot \nabla u dx\right]_S^T
    +  \int_S^T\phi'(t)  \int_{ \Omega }u' h\cdot \nabla u dx dt
    \\
  &\displaystyle +\int_S^T\phi(t)\int_{ \Omega }h\cdot \nabla  u \, f dxdt
  \\
  &\displaystyle -\int_S^T\phi(t)\int_{ \Omega \cap Q_1}\left [\f12 div h(|u'|^2-|\nabla u|^2)
  +\sum_{i,k}\f{\pa u}{\pa x_i}\f{\pa u}{\pa x_k}\f{\pa h_k}{\pa x_i}\right]dxdt\,,\end{array}
  \ee

Using the second multiplier $\f{N-1}{2}\phi(t) u$ for \e{eqn-wave}:
   $$\f{N-1}{2}\int_S^T \phi(t)\int_{\Omega } u(u''-\Delta u -f)dxdt=0\,,$$ 
yields that 
  \be\lab{5.5}
     \begin{array}{rl}
      &\displaystyle \f{N-1}{2}\int_S^T\phi(t)\int_{\Omega } (|\nabla  u|^2-|u'|^2)dxdt
      \\
      =&\displaystyle  -\f{N-1}{2}\left[\phi(t)\int_{\Omega}u' u dx\right]_S^T
      \\
    &\displaystyle+\f{(N-1)}{2}\int_S^T\phi'(t)\int_{\Omega}u' u dx dt
      \\
   &\displaystyle  +\f{N-1}{2}\int_S^T\phi(t)\int_{\Omega} u\,f dxdt\,.
   \end{array}
   \ee
  Adding \e{5.5} to \e{5.4} and using Cauchy-Schwarz and Poincar\'e's inequalities, we obtain for all $\delta_1 >0$ that 
  \be\lab{5.6}\begin{array}{rl}
    &\displaystyle  \int_S^T\phi(t)e(t)dt 
    \\
    =  &\displaystyle\int_S^T\phi(t)\int_{\Omega } \f{|\nabla  u|^2+|u'|^2}{2}dxdt
    \\
     \leq &\displaystyle C  \phi(S)[e(S)+e(T)]+C\int_S^T-\phi'(t)e(t)dt
      \\
    &\displaystyle+\delta_1 \int_S^T\phi(t) e(t)dt 
    +\f{C}{\delta_1}\int_S^T\phi(t)\int_{\Omega}|f |^2 dxdt
    \\
   &\displaystyle +C\int_S^T\phi(t)\int_{\Omega\cap Q_1}(|u'|^2+|\nabla u|^2)dxdt
   \end{array}\,.  \ee

Compared to the desired estimate \e{est-lemma}, the term concerning $|\nabla u|^2$
 on the right hand of \e{5.6} is crucial.  We just follow the techniques developed in \cite{Mart} to deal with this term.

Since $\overline{\mathbb{R}^N\backslash Q_2}\cap\overline{Q}_1=\emptyset$, there
exists a cut-off function $\xi\in C^\infty_0(\mathbb{R})$ such that
 \be \label{xi} 
   0\leq  \xi\leq  1,\quad\xi=1\quad\text{on}\ Q_1,\quad\xi=0 \quad
     \text{on}\ \mathbb{R}^N\backslash Q_2.
     \ee
Applying now the multiplier $\phi(t)\xi(x) u$ to  \e{eqn-wave} gives, after integration by parts, that 
  $$\begin{array}{rl}
    &\displaystyle \int_S^T \phi(t)\int_{\Omega} \xi|\nabla u|^2dxdt
      \\
      =&\displaystyle\int_S^T \phi(t)\int_{\Omega} (\xi|u'|^2+\f12u^2 \Delta \xi)dxdt
      \\
   &+\displaystyle  \int_S^T \phi'(t)\int_{\Omega} \xi uu'dxdt
     -\left[\phi(t)\int_{\Omega}\xi uu'dx\right]_S^T
     \\
    &\displaystyle+\int_S^T \phi(t)\int_{\Omega}\xi u \, fdxdt.
    \end{array}$$
Then it follows from  the definition of $\xi$ that
  \be\lab{5.8}
   \begin{array}{rl}
      \displaystyle\int_S^T\phi(t)\int_{\Omega\cap Q_1}|\nabla u|^2dxdt
   \leq  &\displaystyle \int_S^T \phi(t)\int_{\Omega} \xi|\nabla u|^2dxdt
     \\
     \leq & \displaystyle C \phi(S)[e(S)+e(T)]+C\int_S^T -\phi'(t)e(t)dt
      \\
       & \displaystyle +C\int_S^T\phi(t)\int_{\Omega}|f|^2dxdt
           \\ 
      & \displaystyle+C\int_S^T\phi(t)\int_{\Omega\cap Q_2}(|u'|^2+|u|^2)dxdt\,.
       \end{array}\ee

Now it remains to estimate the term concerning $|u|^2$ in \e{5.8}.
Since $\overline {R^N\backslash \omega}\cap\overline{Q}_2=\emptyset$,
there exists a function $\beta \in C^\infty_0(\mathbb{R})$ such that
   \be \lab{beta}
  0\leq  \beta\leq  1,\quad\beta=1\quad\text{on}\ Q_2,\quad\beta=0 \ \  \text{on}\ \mathbb{R}^N\backslash \omega \,.
   \ee
Fix the $t$ variable and consider the solution $z$ of the following elliptic
problem in space:
   $$\left\{\begin{array}{l} 
   \Delta z=\beta(x)u,\quad\quad  \text{in}\  \Omega,
   \\
  z=0,  \quad\quad\quad\text{on}\ \Gamma.\end{array}\right.
  $$ 
  Hence, $z$ and $z'$ satisfy the following estimates
      \begin{align}\lab{5.10}
         & \|z\|_{  L^2 }\leq  C \displaystyle\int_\Omega \beta(x)|u|^2dx.
         \\     \lab{5.11}  
         & \|z'\|^2_{  L^2 }\leq  C\displaystyle\int_\Omega \beta(x)|u'|^2dx.
\end{align}

Applying  the multiplier $\phi(t) z$  to \e{eqn-wave} gives, after integration by parts, that 
    $$\begin{array}{rl}
     &\displaystyle\int_S^T\phi(t)\int_{\Omega}  \beta(x) |u|^2 dxdt
     \\
     =&\displaystyle\left[ \phi(t)\int_{\Omega} zu'dx\right]^T_S -\int_S^T \phi' (t)\int_{\Omega} zu' dxdt
     \\
     &\displaystyle+ \int_S^T\phi(t)\int_{\Omega}(-z'u'-z f) dxdt\,.
     \end{array}$$
  Hence, using the estimates  \e{5.10}-\e{5.11} in the above relation, and noting the definition of $\beta$,  we obtain for all $\delta_2>0$
  \be\lab{5.12}
   \begin{array}{rl}
      \displaystyle\int_S^T\phi(t)\int_{\Omega\cap Q_2}  |u|^2 dxdt   
    \leq      &\displaystyle\int_S^T\phi(t)\int_{\Omega}  \beta(x) |u|^2 dxdt
       \\
     \leq & \displaystyle C \phi(S)[e(S)+e(T)]+C\int_S^T -\phi' (t)e(t)dt
      \\ 
      &\displaystyle+\delta_2\int_S^T\phi(t)e(t)dt+\f{C}{\delta_2}\int_S^T \phi(t)\int_{\omega}|u'|^2dxdt
      \\ 
      & \displaystyle +\f{C}{\delta_2}\int_S^T\phi(t)\int_{\Omega}|f|^2dxdt\,.
      \end{array}\ee
Inserting \e{5.8} and \e{5.12} in \e{5.6}, then  choosing finally $\delta_1$ and $\delta_2$ sufficiently small, 
we obtain \e{est-lemma} which ends the proof of Lemma \ref{lem-wave}.      \qed

\section*{Acknowledgments}
\small
The authors are thankful to  the support  of the ERC advanced grant 266907 (CPDENL)
and the hospitality of the Laboratoire Jacques-Louis Lions of Universit\'{e} Pierre et Marie Curie. 
Fatiha Alabau-Boussouira was supported by the LIASFMA for a visit at Fudan University during August 2014.
Zhiqiang Wang was partially supported by the National Science Foundation of China (No. 11271082)
and by the State Key Program of National Natural Science Foundation of China (No. 11331004).

%
%

\end{document}